\newif\ifarxiv
\arxivtrue

\ifarxiv
  \documentclass{scrartcl}
\usepackage{algorithm,natbib,amsthm,pdfpages}

\newtheorem{theorem}{Theorem}[section]
\newtheorem{proposition}[theorem]{Proposition}

\newtheorem{corollary}[theorem]{Corollary}

\newcommand{\email}[1]{\protect\href{mailto:#1}{#1}}

\newcommand{\buildtitle}{
\title{Probabilistic Gradients for Fast Calibration of Differential Equation Models.
\thanks{\textbf{Funding} JC was supported by Wave 1 of The UKRI Strategic Priorities Fund under the EPSRC Grant EP/T001569/1, particularly the ``Digital Twins for Complex Engineering Systems'' theme within that grant, and The Alan Turing Institute. ABD was supported by the Lloyds Register Foundation Programme
on Data Centric Engineering and by The Alan Turing Institute under the EPSRC grant
[EP/N510129/1].}}

\author{Jon Cockayne\thanks{The Alan Turing Institute (\email{jcockayne@turing.ac.uk})} \and Andrew B. Duncan\thanks{Imperial College London and The Alan Turing Institute (\email{a.duncan@imperial.ac.uk})}}
}

\newcommand{\thebibstyle}{abbrvnat}
\newcommand{\posttitle}{}
\else
  \documentclass[review]{siamart190516}
  \input{defs_siam}
\fi

\usepackage{amsmath,amssymb,bm,mathtools,color,subcaption}
\usepackage{thmtools}
\usepackage{xr-hyper}
\usepackage{hyperref}
\usepackage[capitalize]{cleveref}
\usepackage{algorithm}
\usepackage{algpseudocode}

\newcommand{\set}[1]{\{ #1 \}}
\newcommand{\Set}[1]{\left\{ #1 \right\}}

\newcommand{\reals}{\mathbb{R}}
\newcommand{\wrt}{\textrm{d}}
\newcommand{\abs}[1]{|#1|}
\newcommand{\Norm}[1]{\left\|#1\right\|}
\newcommand{\norm}[1]{\|#1\|}

\newcommand{\inner}[1]{\langle #1 \rangle}
\newcommand{\Inner}[1]{\left\langle #1 \right\rangle}
\newcommand{\quark}{\setbox0\hbox{$x$}\hbox to\wd0{\hss$\cdot$\hss}}

\newcommand{\defeq}{\coloneqq}

\newcommand\numberthis{\addtocounter{equation}{1}\tag{\theequation}}

\newcommand{\COMMENT}[1]{}

\newcommand{\adj}[1]{#1^\dagger}
\renewcommand\thmcontinues[1]{Continued}

\ifpdf
\hypersetup{
  pdftitle={Probabilistic Gradients for Fast Calibration of Differential Equation Models},
  pdfauthor={J. Cockayne and A. B. Duncan}
}
\fi

\newtheorem{example}[theorem]{Example}

\buildtitle

\begin{document}

\maketitle

\begin{abstract}
Calibration of large-scale differential equation models to observational or experimental data is a widespread challenge throughout applied sciences and engineering.   A crucial bottleneck in state-of-the art calibration methods is the calculation of local sensitivities, i.e.~derivatives of the loss function with respect to the estimated parameters, which often necessitates several numerical solves of the underlying system of partial or ordinary differential equations.  In this paper we present a new probabilistic approach to computing local sensitivities.  The proposed method has several advantages over classical methods.  Firstly, it operates within a constrained computational budget and provides a probabilistic quantification of uncertainty incurred in the sensitivities from this constraint.  Secondly, information from previous sensitivity estimates can be recycled in subsequent computations, reducing the overall computational effort for iterative gradient-based calibration methods.  The methodology presented is applied to two challenging test problems and compared against classical methods. 

\end{abstract}

\posttitle

\section{Introduction}

Complex systems arising in applied sciences and engineering are often modelled by systems of coupled ordinary or partial differential equations (ODEs or PDEs) derived from the underlying physical principles.   
Typically, the specific model behaviour will depend on a vector of parameters which must be \emph{calibrated} to observations of system.  
A major challenge in calibration is the high computational cost associated with numerically solving the mathematical model for a given value of the parameters.  
This is particularly relevant for large-scale models incorporating multi-physics and multiscale behaviour, as arise in the context of digital twins \citep{shafto2012modeling}.  
This high cost often precludes the use of many iterative methods for calibration, including both optimisation methods and Bayesian approaches that use sampling algorithms such as Markov chain Monte-Carlo (MCMC).
Each of these requires at least one solve of the governing equations per iteration of the algorithm.  
In practice MCMC often requires on the order of $10^5$ model evaluations \citep{geyer2011introduction}.

The calibration of differential equation models to observed data can be formulated as a constrained optimisation problem \citep{biegler2007real,gunzburgerperspective,ito2008lagrange}, which is solved using deterministic or stochastic optimisation methods. 
Most fundamental optimisation methods\footnote{i.e.~any method in the Newton family of optimisation methods.} either require or are accelerated by access to derivatives of the functional to be minimised,  so that the solver for the underlying equations must be augmented with a routine that provides the derivative of the solution with respect to model parameters.  
Employing an approximation of the gradient, such as a finite-difference approximation, may seem attractive due to ease of implementation, but obtaining accurate approximations can be challenging and, when the parameter dimension is large, computationally expensive.
Thus it is usually preferable to obtain derivatives using \emph{first-order sensitivity analysis}, which expresses the derivatives as the solution of an auxiliary system of differential equations known as the \emph{sensitivity equations}.  
While the sensitivity equations are linear, they depend on the solution of the underlying equations and so must typically be solved numerically.
Thus, computing the sensitivities is at least as expensive as solving the system itself.  
Further, sensitivities must be computed for every parameter value at which a gradient evaluation is required, making them prohibitively expensive for use in optimisation methods, where gradients are typically required over a large sequence of parameter values.

In the context of model calibration and uncertainty quantification, Gaussian processes (GPs) are often used as surrogate models for the solution of the underlying equations with the aim of making the calibration of such models tractable \citep{higdon2004combining,stuart2018posterior}.   
This approach is advantageous as derivatives of the GP posterior mean can usually be computed explicitly, permitting the use of gradient-based optimisation methods and sampling.    
While GP surrogate models do provide an effective approach to calibrating black-box computer codes where little is known about the structure of the underlying model, this comes at the price of \emph{data-efficiency}.
Information about the gradient can only be obtained from multiple function evaluations near to the location of the required gradient, so again numerous evaluations may be required, particularly if the dimension of the parameter space is high.

These highlighted issues motivate the novel approach to computing sensitivities for optimisation problems presented in this paper.
Our proposed approach is able to bridge the gap between the classical approach of numerically solving the sensitivity equations and the purely data-driven surrogate model approach.
This is achieved by introducing a nonparametric Gaussian process model for the solution of the sensitivity equations that is defined over the entire parameter space.  
The output is a posterior distribution on the space of vector fields in parameter space, whose mean can be interpreted as an estimate of the local sensitivity across multiple parameter locations and whose variance controls the error in this estimate under regularity assumptions.  

This approach offers various advantages to the state-of-the-art approaches:  Firstly, the computational cost of the method can be carefully controlled by the user, either to attain a desired level of accuracy as measured by the ``width'' of the posterior distribution or to fit within a given fixed computational budget.   
Secondly, estimates of gradients at multiple parameter locations are able to share information between them to provide accurate gradient approximations without necessitating additional numerical solves of the underlying PDE model.  
Thirdly, the posterior distribution can be efficiently updated when a gradient evaluation at a new parameter value is required.  
These three advantages are particularly pertinent to model calibration methods which require multiple gradient evaluations along a trajectory.

Besides the immediate application to calibration of PDE models, the efficient approximation of sensitivities for large scale PDE models is of independent interest, with wide ranging applications including model order reduction \citep{pulch2015sensitivity}, shape optimisation \citep{newman1999overview} and uncertainty quantification \citep{arriola2009sensitivity}.   
The  probability distribution output from our new approach has a rigorous Bayesian  interpretation, allowing it to be composed within inference and computation pipelines in a coherent manner to enable propagation of uncertainty.

\subsection{Related Work}

ODE- or PDE-constrained optimisation problems are a class of control problem in which the cost function involves the solution of a partial differential equation posed on a domain $D \subseteq \mathbb{R}^d$.  
Classically, such optimisation problems arise in the context of design and control of engineering systems, for example in optimal topological design, shape design and optimal control of dynamic systems.
See \citet{Herzog2010} for a review of optimisation algorithms for use in this context.
Further, these problems arise naturally in the context of Bayesian inverse problems and model calibration.  
In particular, variational approaches to data assimilation for weather prediction can be naturally rephrased as PDE-constrained optimisation problems \citep{fisher2009data}.  

Sensitivity analysis seeks to quantify the dependence of a function $g(p)$ on perturbations of the problem data or parameters $p \in P$.   
Broadly speaking, we distinguish between \emph{global} sensitivity analysis, which quantifies how input variability influences output variability of a model, and \emph{local} sensitivity analysis which assesses the influence of infinitesimal input perturbations on model output.  
The former is typically assessed in terms of variance, classically using variants of Sobol$'$ indices \citep{Sobol2001}.  
By contrast, local sensitivity analysis involves the calculation of partial derivatives of function outputs with respect to parameters.  
Local sensitivity analysis plays a fundamental role in the context of ODE- or PDE-constrained optimisation \citep{bonnans2013perturbation}.  
In this setting, let $g(u, p)$ denote the real-valued objective function for the optimisation problem, that depends on the solution $u(p)$ of a differential equation for a given parameter value $p \in P$.
Then we seek to compute the total derivative $\frac{\wrt g}{\wrt p}(p)$, which constitutes the local sensitivities.

Generally speaking there are two approaches to computing such derivatives: the \emph{forward} or \emph{direct} method and the \emph{adjoint} method.  
In the forward method, supposing that $p \subseteq \reals^m$, the underlying equations are differentiated with respect to $p_1,\dots,p_m$ to obtain a system of $m$ equations for the sensitivities.  
The adjoint method originates in the theory of Lagrange multipliers in optimisation, and involves solving an auxiliary adjoint equation for the Lagrange multiplier $\lambda$ from which the sensitivities can be directly computed.
Given that the forward approaches involves solving a system of $m$ equations while the adjoint approach involves solving only a single equation, the latter approach can be significantly more efficient for large $m$ \citep{sengupta2014efficient}.

The computational cost of solving optimisation problems involving large-scale ODE or PDE models has motivated the use of \emph{surrogate models}; approximations of the underlying model that can be evaluated at lower computational cost.   
Proposed approaches include using reduced order modelling based on reduced basis methods or proper orthogonal decompositions \citep{benner2015survey, benner2014model}.  
These surrogate approaches are motivated by the fact that the adjoints, and therefore the gradients, of the low-dimensional surrogate model can obtained efficiently. 
Recent efforts involve combining neural network models with low-dimensional physical models to obtain efficient and accurate surrogate models \citep{drohmann2015romes,san2018extreme, san2018neural,hartman2017deep,sheriffdeen2019accelerating}.  
Again, these methods exploit the fact that gradients of neural network models can be obtained efficiently through back-propagation.

Gaussian processes (GPs) have been widely used to provide black-box emulation of computationally expensive codes \citep{sacks1989design}, with \citep{Kennedy2001} providing a mature Bayesian formulation to the methodology.  
Emulation methods based on GPs are now widespread and find uses in numerous applications ranging from computer code calibration \citep{higdon2004combining}, uncertainty analysis \citep{oakley2002bayesian} and MCMC \citep{lan2016emulation,CLEARY2020109716}.  
Among the first papers to consider application of emulation within sensitivity analysis was \citet{Oakley2004}, which extended the work of \citet{Kennedy2001} to computation of variance-based global sensitivities.
Subsequent work by \citet{Jin2004} considered a similar approach that exploited a tensor-product kernel to simplify the integration problems required, though this work did not consider the posterior covariance in their estimator.  
See \citet{Cheng2020} for a more extensive review of emulation-based global sensitivity analysis techniques, and \citet{Girard2016,Beddows2017,Renardy2018} for a survey of applications of such approaches.
One could envisage an analogous emulation strategy for local sensitivity analysis of computationally expensive models that involves first constructing an emulator  $\hat{g}$ of the objective function $g$ and then evaluating the derivative $\frac{\wrt \hat{g}}{\wrt p}(p^*)$ which, assuming a conducive emulator, can be computed at a lower cost than the derivative of $g$ itself.
A notable disadvantage of this approach is that to approximate local sensitivities in this way  would require \emph{global} information about $g$, since unless a highly structured surrogate model is used little information can be obtained about $\frac{\wrt g}{\wrt p}(p)$ from the single evaluation $g(p)$.

The method proposed in this paper aims to bridge the gap between classical numerical approaches and emulation-based approaches to calculating sensitivities within optimisation problems.   
The proposed method can be interpreted as a Bayesian probabilistic numerical method \citep{Cockayne2019} for the solution of the forward or adjoint sensitivity equations over $D \times P$.
It is similar to the probabilistic meshless methods for solutions of PDEs presented in \citet[Chapter 5]{CockayneThesis}, but extended across parameter space.
This formalism presents several advantages.
Firstly it permits a high level of adaptivity, in that the solution can be refined over both $P$ and $D$ to increase accuracy either globally over parameter space $P$, or locally for particular value of the parameters $p \in P$.
Secondly, subject to regularity assumptions, estimates of the gradient at a parameter $p$ may re-use information from nearby gradient evaluations, exploiting the smoothness of the sensitivity equations to reduce the computational effort required for accurate gradient estimates at $p$ when nearby gradients have already been evaluated.  
Thirdly, gradient estimates can be updated efficiently, allowing the adaptivity and smoothness properties mentioned to be exploited within algorithms that depend upon local sensitivities, such as gradient-based optimisation algorithms.

\subsection{Contributions}

The main contributions of the paper are as follows:
\begin{itemize}
	\item We develop a probabilistic framework for computing gradients for differential equation models.
	\item We study the theoretical properties of this method, in particular its robustness to discretisation error.
	\item We demonstrate how the inferred gradients can be leveraged in optimisation problems. 
	\item The results are demonstrated on a number of model problems to analyse the method's performance in comparison to classical approaches.
\end{itemize}

\subsection{Structure of the Paper}
The paper proceeds as follows.
In \cref{sec:classical_sensitivity_analysis} the classical approach to computing local sensitivities is formulated with examples of application to the problem of computing sensitivities for a simple PDE.
\cref{sec:probabilistic} presents the novel probabilistic approaches and provides theoretical results relating to their accuracy and stability.
\cref{sec:optimisation} discusses the use of the probabilistic methods introduced in optimisation problems, and the empirical performance of these methods is assessed in \cref{sec:applications}.
We conclude with some discussion in \cref{sec:conclusion}.
The supplementary material contains the proofs required for the paper in \cref{sec:proofs}.

\subsection{Notation}

Let $W^{k,p}(D)$ denote the Sobolev space in which each function has $k$ weak derivatives with finite $L^p(D)$ norm.
We will use the notation $H^k(D) = W^{k,2}(D)$.
Further let $H^k_0(D)$ denote the subset of $H^k(D)$ for which all $f \in H^k_0(D)$ have $f = 0$ on $\partial D$ and $H^{-k}(D)$ to be the dual of $H_0^k(D)$.
For two normed spaces $\mathcal{U}$, $\mathcal{V}$ we will use the notation $\mathcal{L}(\mathcal{U}, \mathcal{V})$ to denote the set of all bounded linear operators from $\mathcal{U}$ to $\mathcal{V}$.
For the set of all bounded linear functionals on $\mathcal{U}$ we will use the notation $\mathcal{U}^* = \mathcal{L}(\mathcal{U}, \reals)$.
When $\mathcal{U}$ is a set of functions on some domain $D$ we will use the notation $\delta[x]$ to denote the evaluation functional for the point $x \in D$, i.e.~$\delta[x](u) = u(x)$.

When both $\mathcal{U}$ and $\mathcal{V}$ are Hilbert spaces, for an operator $A \in \mathcal{L}(\mathcal{U}, \mathcal{V})$ let $\adj{A} \in \mathcal{L}(\mathcal{V}, \mathcal{U})$ denote the adjoint of $A$.
For $A \in \mathcal{L}(\mathcal{U}, \mathcal{U})$, recall that the trace of $A$ is defined as $\textup{trace}(A) = \sum_{i=1}^\infty \inner{A e_i, e_i}$
where $(e_i)_{i=1}^\infty$ is an arbitrary orthonormal basis of $\mathcal{U}$. 

Several operator norms will be required.
For an operator $A : \mathcal{U} \to \mathcal{V}$ we will denote the operator norm by $\norm{A}_{\mathcal{U} \to \mathcal{V}} = \sup_{u \in \mathcal{U}} \norm{A u}_{\mathcal{V}} / \norm{u}_{\mathcal{U}}$.
When $\mathcal{U} = \mathcal{V}$ we will simply use the notation $\norm{A}_{\mathcal{U}}$.
The trace norm is given by $\norm{A}_{\textsf{tr}} = \textup{trace}([\adj{A}A]^\frac{1}{2})$ while the Hilbert-Schmidt norm is given by $\norm{A}_{\textsf{HS}} = \textup{trace}(\adj{A} A)^\frac{1}{2}$.
Recall that $\norm{A}_{\mathcal{U}\to\mathcal{V}} \leq \norm{A}_\textsf{HS} \leq \norm{A}_\textsf{tr}$.

\subsubsection{Fr{\'e}chet Derivatives}
Of central importance to the paper is the concept of a Fr{\'e}chet derivative.
Let $\mathcal{U}$ and $\mathcal{V}$ each be normed spaces and consider a function $f : \mathcal{U} \to \mathcal{V}$.
When it exists, Fr{\'e}chet derivative of $f$ at $u \in \mathcal{U}$ is defined to be the operator $\frac{\wrt f}{\wrt u}[u] \in \mathcal{L}(\mathcal{U}, \mathcal{V})$ that satisfies
\begin{equation} \label{eq:frechet_one_var}
	\lim_{\norm{h} \to 0} \frac{\norm{f(u + h) - f(u) - \frac{\wrt f}{\wrt u}[u]h}}{\norm{h}} = 0
\end{equation}
where the notation $\norm{h} \to 0$ is a shorthand for the requirement that the limit exist uniformly across sequences $(h_n)$ in $\mathcal{U}$ such that $\norm{h_n} \to 0$ as $n\to \infty$.
It is important to observe that $\frac{\wrt f}{\wrt u}[u]$ is a linear operator that depends upon $u$, so that  $\frac{\partial f}{\partial u}[u](v)$ is the Fr{\'e}chet derivative at the location $u \in \mathcal{U}$ in the direction $v \in \mathcal{U}$.  For a function $f : \mathcal{U}\times \mathcal{V} \to \mathcal{W}$ the partial Fr{\'e}chet deriviative is defined analagously to be the operator $\frac{\partial f}{\partial u} \in \mathcal{L}(\mathcal{U} \times \mathcal{V}, \mathcal{W})$ that satisfies
\begin{equation*}
	\lim_{\norm{h} \to 0} \frac
		{\norm{f(u + h, v) - f(u, v) - \frac{\partial f}{\partial u}[u, v](h)}}
		{\norm{h}} = 0
\end{equation*}
whenever the above limit exists.  Finally, consider the case where the function $u$ depends on $v$. Let $f : \mathcal{U} \times \mathcal{V} \to \mathcal{W}$, suppose that $\mathcal{U}$ is a space of functions with domain $\mathcal{V}$.
Then, when it exists, the Fr{\'e}chet derivative of $f$ with-respect-to $v$ is the operator $\frac{\wrt f}{\wrt v} \in \mathcal{L}(\mathcal{U} \times \mathcal{V}, \mathcal{W})$ that satisfies
\begin{equation} \label{eq:frechet_derivative_total}
	\lim_{\norm{h} \to 0} \frac
		{\norm{f(u(v+h), v+h) - f(u, v) - \frac{\wrt f}{\wrt v}[u(v), v](h)}}
		{\norm{h}} = 0.
\end{equation}
This will sometimes be referred to as the \emph{total} Fr{\'e}chet derivative of $f$.

\section{Background} \label{sec:classical_sensitivity_analysis}

In this section a formal presentation of local sensitivity analysis is provided.
In \cref{sec:background:problem} the problem is introduced, while \cref{sec:background:forward,sec:background:adjoint} present forward and adjoint sensitivity analysis, respectively.
Lastly in \cref{sec:background:pn_pdes} we will briefly discuss probabilistic numerical methods for the solution of PDEs, and discuss their similarity to this work.

\subsection{Local Sensitivity Analysis} \label{sec:background:problem}

We begin by introducing the relevant spaces for the problem.
Let $\mathcal{U}$, $P$, $\mathcal{F}$ and $\mathcal{G}$ each be real-valued Banach spaces.
In this paper it will be assumed that $\mathcal{U}$ and $\mathcal{F}$ are infinite-dimensional spaces of functions defined on spatial domain $D$, with $\mathcal{U}$ referred to as the \emph{solution space} and $\mathcal{F}$ as the \emph{constraint space}.
Define  $\mathcal{U}_P$ to be a space of real-valued functions on $D \times P$ with the property that $u(\quark, p) \in \mathcal{U}$ for all $p \in P$, and let
$\mathcal{U}_{\partial P} = \Set{\frac{\partial u}{\partial p} : u \in \mathcal{U}_P}.$
The set $\mathcal{F}_P$ is defined analogously.
The \emph{parameter space} $P$ may be finite- or infinite-dimensional.
The space $\mathcal{G}$ will be referred to as the \emph{quantity of interest (QoI) space}, and will be assumed to be finite-dimensional.
In particular it will often be the case that $\textup{dim}(\mathcal{G}) = 1$, though we note that this is not required for the presentation below.

Two functions define the problem.
The function $F : \mathcal{U}\times P \to \mathcal{F}$ is referred to as the \emph{constraint function}, and loosely speaking this encapsulates all of the constraints that must be satisfied in order for a pair $(u, p) \in \mathcal{U} \times P$ to constitute a solution to the PDE.
The function $g : \mathcal{U} \times P \to \mathcal{G}$ is referred to as the \emph{QoI function}, and this describes a typically low-dimensional quantity of interest derived from the solution; in the context of optimisation problems this will generally be the objective function whose minimiser is sought.

More formally, $F$ is such that for each $p \in P$ there is a unique $u^\dagger \in \mathcal{U}_P$ that satisfies
$F(u^\dagger(\quark, p), p) = 0$
for each $p \in P$.
For convenience, define the \emph{parameter-to-solution map} $U : P \to \mathcal{U}$ which provides the solution to the underlying differential equation for a particular value of the parameter, i.e.~$U(p) = u^\dagger(\quark, p)$.
As a result, the equation $F(U(p), p) = 0$ is automatically satisfied for all $p \in P$.

It will be assumed the partial Fr\'echet derivatives of $F$ and $g$ with-respect-to both $u$ and $p$ exist and are tractably computable for all pairs $(u, p) \in \mathcal{U}\times P$.
It will also be assumed that the derivative of $U$ with-respect-to $p$ exists but is not tractable.
Note that this implies the existence of the total derivatives $\frac{\wrt F}{\wrt p}$ and $\frac{\wrt g}{\wrt p}$.
Lastly we assume that $\frac{\partial F}{\partial u}[U(p), p]$ is nonsingular for each $p \in P$.

The objective is to estimate the value of the Fr\'echet derivative
\begin{equation*}
	\frac{\wrt g}{\wrt p}[U(p), p] \in \mathcal{L}(P, \mathcal{G})
\end{equation*}
for a pair $(U(p), p)$.
Note that since the location at which the derivative is taken is $U(p)$, this should be interpreted as a total Fr\'echet derivative in the form of \cref{eq:frechet_derivative_total}.
To fix ideas we consider the following simple parameter sensitivity problem.

\begin{example}[Partial Differential Equation] \label{example:pde}
Let $P \subseteq \reals^n$ be an open set. Consider the following parametrised steady state conductivity model:
\begin{align*}
	-\nabla \cdot (\kappa(x; p) \nabla u(x)) &= f(x) &x&\in D \\
	u(x) &= 0 & x&\in \partial D
\end{align*}
where $f \in H^{-1}(D)$ and $\kappa : D\times P \to \reals^{d\times d}$ satisfies $\lambda_p |e|^2 \leq e\cdot \kappa(x, p)e  \leq \Lambda_p|e|^2$ for all $x \in D$ and $e \in \mathbb{R}^d$ for some $0 < \lambda_p < \Lambda_p < \infty$ for all $p \in P$.  Standard existence theory for elliptic PDEs \citep[Section 6.2, Theorem 3]{Evans2010} states that a weak solution $u \in H^1_0(D)$ exists for every $p \in P$.  For convenience we will suppose that the boundary conditions are implicitly satisfied, i.e.~$\mathcal{U} = H^1_0(D)$. The constraint equation is given by
$F(u, p) = -\nabla \cdot (\kappa(x; p) \nabla u(x)) - f(x)$
so that $\mathcal{F} = H^{-1}(D)$. Suppose that the quantity-of-interest is $g(x) = \norm{u}_2 = \left(\int_D u^2(x) \,\wrt x\right)^\frac{1}{2}$.

\end{example}

Both forward and adjoint sensitivities are computed by first observing that the total derivative of interest, $\frac{\wrt g}{\wrt p}$ satisfies the following identity:
\begin{equation} \label{eq:total_deriv_expansion}
	\frac{\wrt g}{\wrt p}[U(p), p] = \frac{\partial g}{\partial u}[U(p), p]\frac{\wrt U}{\wrt p}[p] + \frac{\partial g}{\partial p}[U(p), p].
\end{equation}
Since it is assumed that $\frac{\partial g}{\partial u}$ and $\frac{\partial g}{\partial p}$ are each analytically tractable, the only remaining quantity that must be computed is $\frac{\wrt U}{\wrt p}$.
The challenge is that since the parameter-to-solution map $U(p)$ is typically inaccessible and must be approximated independently for each $p \in P$, $\frac{\wrt U}{\wrt p}$ is also difficult to compute.
The forward and adjoint approaches handle this intractability in different ways, which will now be presented.

\subsection{Forward Sensitivity Analysis} \label{sec:background:forward}

In forward sensitivity analysis we seek to calculate $\frac{\wrt U}{\wrt p}$ directly.
Note that we have
\begin{equation*}
	\frac{\partial F}{\partial p}[U(p), p] = \frac{\partial F}{\partial u}[U(p), p]\frac{\wrt U}{\wrt p}(p) + \frac{\partial F}{\partial p}[U(p), p], \quad p \in P.
\end{equation*}
Further, since by construction $F(U(p),p) = 0$, we also have that $\frac{\wrt F}{\partial p}[U(p), p] = 0$.
This gives the forward sensitivity equation
\begin{equation} \label{eq:F_deriv_0}
	\frac{\partial F}{\partial u}[U(p), p] \frac{\wrt U}{\wrt p}[U(p), p] = -\frac{\partial F}{\partial p}[U(p), p],\quad p \in P
\end{equation}
which is a linear system whose solution can be computed to determined $\frac{\wrt U}{\wrt p}$, since $\frac{\partial F}{\partial u}$ is assumed to be invertible.
This solution can then be substituted into \cref{eq:total_deriv_expansion} to compute $\frac{\wrt g}{\wrt p}$.

Note that both the operator $\frac{\partial F}{\partial u}[U(p), p]$ and the right-hand-side $-\frac{\partial F}{\partial p}[U(p), p]$ depend both on the parameter value $p$ and the solution $U(p)$.
This has two important consequences.
Firstly, if sensitivities are required at another point $q\neq p$ then the solution $U(q)$ must be recomputed and the forward sensitivity equation \cref{eq:F_deriv_0} must be solved anew to determine $\frac{\wrt U}{\wrt p}[U(q), q]$.
Secondly, for most problems of interest $U(p)$ will not be available explicitly and one must substitute an approximate solution $\hat{U}(p) \approx U(p)$.
This may induce further numerical error, the impact of which must in turn be analysed, but also means that even though \cref{eq:F_deriv_0} is linear, its solution is unlikely to be available in closed-form owing to its dependence on $\hat{U}(p)$.
We now consider the computation of the forward sensitivities for \cref{example:pde}.

\begin{example}[Elliptic PDE: Forward Sensitivity Analysis] \label{example:pde_forward}
We begin by deriving $\frac{\partial F}{\partial p}$. Assume that $\kappa$ is once-differentiable in each coordinate of $p$ and that $\sup_{x\in D}|\partial_{p_i}\kappa(x; p)| < \infty$.   The Frech{\'e}t derivative of $F$ with respect to $p$ at $(U(p), p)$ is defined by
\begin{equation}
\label{eq:elliptic_pde_F_by_p}
	\frac{\partial F}{\partial p}[U(p), p]q %&= -\nabla_p \nabla \cdot \left(\kappa(x; p) \nabla u(x) \right) \\
	= -\sum_{i=1}^{m}
		\nabla \cdot \left(\frac{\partial \kappa}{\partial p_i}(x; p) \nabla U(p)(x) \right) q_i, \quad q \in P.
\end{equation}
From energy estimates for weak solutions of elliptic PDEs, $\nabla U(p)(x) \in L^2(D)$.
For illustration, it is straightforward to show that the RHS of \cref{eq:elliptic_pde_F_by_p} lies in $H^{-1}(D)$.  
The derivative $\frac{\partial F}{\partial u}$ is given by 
\begin{equation} \label{eq:example_df_du}
	\frac{\partial F}{\partial u}[U(p), p](v) = -\nabla \cdot \left(\kappa(x; p) \nabla v(x) \right) \quad v \in \mathcal{U}.
\end{equation}
so that clearly $\frac{\partial F}{\partial u}[U(p), p] \in \mathcal{F}$ since in this case, owing to the linearity of the PDE operator, $\frac{\wrt F}{\wrt u}$ is identical to this operator and independent of both $p$ and $U(p)$, though for general nonlinear problems this will not be the case.
The sensitivities of $U$ with respect to the $p_i$ are therefore defined by the following system of PDEs
\begin{equation} \label{eq:fsa_system}
		-\nabla \cdot \left(\kappa(x; p) \nabla \frac{\wrt U(p)}{\wrt p_i}(x) \right)
		= -\nabla \cdot \left( \frac{\partial \kappa}{\partial p_i}(x; p) \nabla U(p)(x)\right), \quad (x,p) \in D\times P.
\end{equation}
For fixed $p$ system of equations is well-posed, guaranteeing the existence of unique solutions $\frac{\wrt U}{\wrt p_i} \in H_0^1(D)$, $i=1,\dots,m$. 

Once these $m$ PDEs have been solved, the computed solutions can be substituted into $\cref{eq:total_deriv_expansion}$ to determine $\frac{\wrt g}{\wrt p}$.
To accomplish this we are required to compute the derivatives $\frac{\partial g}{\partial u}$ and $\frac{\partial g}{\partial p}$.
Note that in this case $g$ is independent of $p$, and it is further straightforward to show that
\begin{equation*}
	\frac{\partial g}{\partial u}[u](v) = \frac{\wrt g}{\wrt u}[u](v) = \frac{\inner{u, v}_2}{\norm{u}_2}.
\end{equation*}
Once again, note that this is a linear operator in $v$, but is nonlinear in $u$.
We therefore have that
\begin{equation*}
	\frac{\wrt g}{\wrt p_i}[U(p), p] = \frac{1}{\norm{U(p)}_2} \Inner{U(p), \frac{\wrt U}{\wrt p_i}}_2
\end{equation*}
for the derivatives $\frac{\wrt U}{\wrt p_i}$ identified by solution of \cref{eq:fsa_system}.
\end{example}

The central challenge with the forward approach, which motivates the adjoint approach that will be presented in the next section, is the dependence of the forward sensitivity equation \cref{eq:F_deriv_0} on the dimension of the parameter space: solving for $\frac{\wrt U}{\wrt p}$ requires the solution of $\mbox{dim}(P)$ PDEs.  In many practical problems the parameter space is extremely large; thus, a method for computing the sensitivities that is independent of the dimension of the parameter space is also of interest.

\subsection{Adjoint Sensitivity Analysis} \label{sec:background:adjoint}

Adjoint sensitivity analysis begins by introducing the operator $\lambda \in \mathcal{L}(\mathcal{F}, \mathcal{G})$. 
Supposing that $\textup{dim}(\mathcal{G}) = n$, we can express $g$ as $(g_1, \ldots, g_n)$ and consequently $\lambda = (\lambda_1, \ldots, \lambda_n)$ where $\lambda_i \in \mathcal{F}^*$ for $i =1,\ldots, n$.
For fixed $p$, the auxiliary term $\lambda$ is selected to solve
\begin{equation} \label{eq:lambda_system}
	\lambda_i \frac{\partial F}{\partial u}[U(p), p] = \frac{\partial g_i}{\partial u}[U(p), p], \quad i=1,\ldots, n.
\end{equation}
Assuming this is a unique solution $\lambda$ exists, one can then recover the sensitivity of the quantity of interest $g$ as follows
\begin{equation} \label{eq:adjoint}
	\frac{\wrt g_i}{\wrt p} = - \lambda_i \frac{\partial F}{\partial p} + \frac{\partial g}{\partial p},  \quad i=1,\ldots, n.
\end{equation}
which provides a computable expression for the local sensitivities.

We note that compared to \cref{sec:background:forward} which, in the finite-dimensional case, necessitates $m = \mbox{dim}(P)$ solutions of the forward sensitivity equation, the adjoint system requires $n = \mbox{dim}(G)$ solutions of the adjoint sensitivity equation.   
In typical situations where $n \ll m$ then there is a clear computational benefit to this approach.

\begin{example}[Elliptic PDE: Adjoint Sensitivity Analysis] \label{example:pde_adjoint}
Recalling $\frac{\partial F}{\partial u}$ and $\frac{\partial g}{\partial u}$ as derived in \cref{example:pde_forward}, the problem that must be solved to identify $\lambda \in H_0^1(D)$ such that
\begin{equation} \label{eq:asa_lambda}
	\nabla \cdot\left(\kappa^{\top}(x; p) \nabla \lambda(x) \right) = \frac{U(p)}{\norm{U(p)}_2}.
\end{equation}

Once $\lambda$ has been determined, referring again to the derivation in \cref{example:pde_forward} we have that
\begin{equation*}
	\frac{\partial g}{\partial p_i} = -\int \nabla\lambda(x)  \cdot
		\frac{\partial \kappa}{\partial p_i}(x; p) \nabla U(p)(x)\,dx
\end{equation*}
which is real-valued, as required.
Again note that in the equation that determines $\lambda$, $U(p)$ appears on the right-hand-side, so for each value of $p$ for which sensitivities are required the PDE must be solved.
Nevertheless the fact that in this example only a single system needs to be solved for each $p$ makes the adjoint method significantly cheaper to apply when $m= \mbox{dim}(P)$ is large.
\end{example}

In the next section we will describe the new probabilistic approaches to both forward and adjoint sensitivity analysis, each of which operates with a constrained computational budget.

\subsection{Probabilistic Numerical Methods for PDEs} \label{sec:background:pn_pdes}

When applied to PDEs, there is a marked similarity between this work and probabilistic numerical methods\footnote{See \citep{Hennig2015a} for a high-level introduction, and \citep{Oates2019} for a thorough literature review.} applied to linear PDEs.
In this section we will discuss these methods, and the similarity to the present approach.
Broadly speaking these methods begin by placing a Gaussian prior on the function space occupied by the solution to the PDE.
Finite-dimensional information about the unknown solution is then produced by projecting the linear PDE through a set of $d$ functionals, referred to as \emph{information functionals} in this work.
The conjugacy of Gaussian distributions with linear projections can then be exploited to write down the posterior distribution in closed-form.
For a detailed introduction to this perspective see \citep[Chapter 5]{CockayneThesis}, in which it is referred to as the \emph{probabilistic meshless method} (PMM).

% Radial basis function collocation
This approach is equivalent to symmetric collocation with radial basis functions \citep{Wendland2004,Cialenco2012}, in that it is possible to construct the prior such that the posterior mean from PMM coincides with the estimator for the solution of the PDE produced in symmetric collocation.
To our knowledge this approach was first presented in \citep[Chapter 16]{Wendland2004}, and extended in \citep{Cialenco2012} to refine the error analysis, as well as explore applications in stochastic PDEs.
In symmetric collocation the posterior distribution itself is not of interest, but the error analysis that appears in those works is relevant here as it provides an important interpretation for the posterior covariance.
Specifically, the bound that appears in \citep{Wendland2004} connects the error to an object referred to as the \emph{power function}, which can be shown to be directly connected to the posterior covariance that appears in the PMM.

% FEM-based methods
In addition to the PMM, other works that could be interpreted as probabilistic numerical methods for PDEs include a series of papers that introduced \emph{gamblets} for the solution of PDEs with rough coefficients \citep{Owhadi2015,Owhadi2017b,Owhadi2017a}. 
These papers construct a probabilistic solution to the PDE in a broadly similar way to \citep{CockayneThesis}, but with several distinct differences.
Firstly, the probability model is motivated by a game theoretic argument rather than Bayesian reasoning, though the ultimate conditioning procedure arrived at is equivalent.
Secondly, the information about the solution is constructed in a distinctly different way, by projecting the defining equations of the PDE against a hierarchical basis formed by a nested partitioning of the domain, whereas in the PMM and in symmetric collocation it is obtained by evaluating those equations at a set of points referred to as \emph{collocation points}.
However this results in a very different error analysis, since collocation methods typically bound the estimation error in terms of the fill distance of these collocation points, whereas in gamblet-based methods, since there is no analogue of these points, a different approach must be adopted.

% What is the difference
The chief similarities of these approaches to the approach presented in this paper is that, when the system defined by $F$ is a PDE, the sensitivity equations will involve solving a system of PDEs.
In this setting the approach that we describe is similar in principal to the approaches we describe above, in that for a particular choice of prior and information, the method we employ will be equivalent to these methods.
There are several distinct differences however.
Firstly, it is possible that the system described by $F$ is \emph{not} a PDE, and indeed in this work we will explore sensitivity analysis for ODEs in addition to PDEs.
While there exist probabilistic numerical methods for solving ODEs, they typically make approximations to account for nonlinearity which are not required in this work, as the systems which must be solved in sensitivity analysis are linear.
Secondly, in the PDE case we do note make specific assumptions on the form of the information functionals, as these will typically be problem specific.
Thirdly, by formulating the sensitivity equations as a single (degenerate) PDE on the joint space $D\times P$, the continuity of the sensitivities with respect to $p$ is exploited to permit implicit interpolation of the sensitivities across different values of $p$.
And lastly, the focus of this paper is on computing sensitivities, \emph{not} on the solution of the PDE itself, which is assumed to be obtained by some classical numerical solver.

\section{Probabilistic Approaches} \label{sec:probabilistic}

In this section we will present two probabilistic approaches to computing parameter sensitivities.
Each allows a user to restrict the amount of computational effort expended and still obtain an estimate of the sensitivities, while also providing an estimate of the error incurred as a result.
Familiarity with Gaussian processes is assumed for this section; we refer the unfamiliar reader to the introduction given in \citet{Rasmussen2005}; see also \citet{Bogachev:1998} for a more mathematical treatment.

We will assume that there exist reproducing kernel Hilbert spaces (RKHSs) $\mathcal{U}_P'$, $\mathcal{F}_P'$ such that $\mathcal{U}_P'$ is dense in $\mathcal{U}_P$ and $\mathcal{F}_P'$ is dense in $\mathcal{F}_P$.
Let 
$\mathcal{U}_{\partial P} = \Set{\frac{\partial u}{\partial p} : u \in \mathcal{U}_P}$
and let $\mathcal{U}_{\partial P}'$ be defined analogously for $\mathcal{U}_P'$.
It will also be assumed that $g$ is a functional, so that $\mathcal{G} = \reals$; this last assumption can readily be generalised, and is made to simplify the presentation.

\subsection{Probabilistic Forward Sensitivity Analysis} \label{sec:probabilistic:forward}

We first consider forward sensitivity analysis.
We begin by modelling prior uncertainty about $\frac{\partial U}{\partial p}$ with the random variable $X_F$, distributed as $X_F \sim \mu_F = \mathcal{N}(a_F, C_F)$, where $a_F \in \mathcal{U}_{\partial P}'$ and $C_F :\mathcal{U}_{\partial P}' \to \mathcal{U}_{\partial P}'$ is a positive-definite covariance operator.
It will be assumed that $\mu_F(\mathcal{U}_{\partial P}) = 1$.
When $\textup{dim}(P) < \infty$ this prior takes the form of a vector-valued Gaussian process prior \citep{Alvarez2012}.
In the infinite-dimensional setting, we note that a discretisation of the parameter space will nevertheless be required for computational purposes, resulting in a parameter space that is effectively finite-dimensional, though a finite-dimensional parameter space is not strictly required for the theoretical results presented herein.

To obtain a posterior belief over the forward sensitivities, this prior will be conditioned on observations of \cref{eq:F_deriv_0}.
Let $\tilde{\mathcal{I}}_{F,1},\dots,\tilde{\mathcal{I}}_{F,d}$ be a such that $\tilde{\mathcal{I}}_{F,j} \in (\mathcal{F}^m)^*$ for $j=1,\dots,d$ and let $\set{p_1, \dots, p_d} \subset P$.
Let $X_F$ be a random variable with law $\mu_F$.
Note that the prior distribution $\mu_F$ implies a prior distribution over $\frac{\wrt g}{\wrt p}$ by projecting through the linear map given in \cref{eq:total_deriv_expansion}; this will be denoted $\nu_F$.
By applying each operator $\tilde{\mathcal{I}}_{F,j}$ to \cref{eq:F_deriv_0} we obtain
\begin{equation} \label{eq:forward_info}
	\tilde{\mathcal{I}}_{F,j} \frac{\partial F}{\partial u}[U(p_j), p_j] X_F = -\tilde{\mathcal{I}}_{F,j} \frac{\partial F}{\partial p}[U(p_j), p_j]
\end{equation}
which, under the assumptions made at the start of this section, yields the information
\begin{equation*}
	f_{F,j} = -\tilde{\mathcal{I}}_{F,j} \frac{\partial F}{\partial p}[U(p_j), p_j]
\end{equation*}
where $f_{F,j} \in \reals$.
Let $f_F \in \reals^d$ be the vector with $[f_F]_j = f_{F, j}$. 

It is more mathematically convenient to think of the $\tilde{\mathcal{I}}_{F,j}$ in terms of functionals defined on $\mathcal{U}_{\partial P}$.
To this end, let $\mathcal{I}_{F, j} \in \mathcal{U}_{\partial P}^*$ be defined by 
\begin{equation*}
	\mathcal{I}_{F, j} \frac{\partial u}{\partial p} = \tilde{\mathcal{I}}_{F,j} \frac{\partial F}{\partial u}[U(p_j), p_j] \frac{\partial u}{\partial p}(\quark, p_j) .
\end{equation*}
We refer to $\mathcal{I}_{F,1},\dots,\mathcal{I}_{F,d}$ as the \emph{information functionals}, and will assume that the information functionals are linearly independent.

The posterior is obtained by conditioning the prior on the information functionals.
First, introduce the operator $\mathcal{I}_F : \mathcal{U}_P \to \reals^d$, given by
\begin{equation*}
	\mathcal{I}_F \frac{\partial u}{\partial p} = \begin{bmatrix}
		\mathcal{I}_{F,1} \frac{\partial u}{\partial p} \\
		\vdots\\
		\mathcal{I}_{F,d} \frac{\partial u}{\partial p}
	\end{bmatrix}.
\end{equation*}
Then we seek to compute $X_F | \mathcal{I}_F X_F = f_F$.
Owing to the linearity of $\mathcal{I}_F$, the resulting posterior distribution is again Gaussian and is given in the following proposition.

\begin{proposition}[Probabilistic Forward Sensitivity Analysis] \label{prop:forward_posterior}
	The posterior $X_F | \mathcal{I}_F X_F = f_F$ has law $\bar{\mu}_F$ given by
	\begin{align*}	
		\bar{\mu}_F &= \mathcal{N}(\bar{a}_F, \bar{C}_F) \\
		\bar{a}_F &= a_F + C_F\adj{\mathcal{I}_F}[\mathcal{I}_FC_F \adj{\mathcal{I}_F}]^{-1} (f_F - \mathcal{I}_F a_F ) \\
		\bar{C}_F &= C_F - C_F\adj{\mathcal{I}_F}[\mathcal{I}_FC_F \adj{\mathcal{I}_F}]^{-1}\mathcal{I}_FC_F .
	\end{align*}
	The implied posterior distribution over $\frac{\wrt g}{\wrt p}$, denoted $\bar{\nu}_F$, is given by
	\begin{align*}
		\bar{\nu}_F &= \mathcal{N}(\bar{g}_F, \bar{G}_F) \\
		\bar{g}_F(p) &= \frac{\partial g}{\partial u}[U(p), p](\bar{a}_F(\quark, p)) + \frac{\partial g}{\partial p}[U(p), p] \\
		\bar{G}_F(p,p') &= \frac{\partial g}{\partial u}[U(p), p] \delta[\quark, p] \bar{C}_F \adj{\delta[\quark, p']} \adj{\frac{\partial g}{\partial u}[U(p'), p']}.
	\end{align*}

\end{proposition}

An important note is that even when underlying system described by $F$ is nonlinear, the posterior distribution remains Gaussian owing to the linearity of the Fr{\'e}chet derivatives.
Choice of prior mean and covariance is highly problem specific, and will be discussed for the specific examples considered in this paper in \cref{sec:applications}.
Next we turn to the adjoint approach.

\subsection{Probabilistic Adjoint Sensitivity Analysis} \label{sec:probabilistic:adjoint}

For the adjoint problem, the system that must be solved is now \cref{eq:lambda_system}.
Since $\mathcal{F}_P'$ is assumed to be an RKHS, due to the representer theorem \citep[see e.g.][Section 4.4]{Berlinet2004} we have $\lambda f = \Inner{f, \beta}_\mathcal{F}$, where $f,\beta \in \mathcal{F}_P'$, 

The proposed approach is as in the previous section.
We model uncertainty in $\beta$ with the random variable $X_A$, whose law is $\mu_A = \mathcal{N}(a_A, C_A)$,
where $a_A \in \mathcal{F}_P'$ and $C_A : \mathcal{F}_P'\to \mathcal{F}_P'$ is a positive-definite covariance operator.
Note that this again implies a distribution $\nu_A$ over $\frac{\wrt g}{\wrt p}$ by projecting through the linear map
\begin{equation*}
	\mathcal{J}[p](\beta) = \Inner{\frac{\partial F}{\partial p}[U(p), p], \beta(\quark, p)}_\mathcal{F}.
\end{equation*}
An important remark, however, is that unless $\mu_A$ and $\mu_F$ are chosen carefully, the implied distributions $\nu_A$ and $\nu_F$ will not be equal.

To define the information functionals let $\set{(e_1, p_1),\dots,(e_d, p_d)} \subset \mathcal{U} \times P$. 
Then
\begin{equation} \label{eq:information_adjoint}
	\mathcal{I}_{A,j} \beta = \Inner{\frac{\partial F}{\partial u}[U(p_j), p_j](e_j), \beta(\quark; p_j)}_{\mathcal{F}},
\end{equation}
so that $\mathcal{I}_{A,j} \in \mathcal{F}^*_P$.
Furthermore note that the information $f_{A,j} \defeq \frac{\partial g}{\partial u}[U(p_j), p_j](e_j)$ is clearly computable.
Let $f_A$ and $\mathcal{I}_A$ be defined analogously to previous sections; then the posterior on $\beta$ is given in the following proposition.

\begin{proposition}[Probabilistic Adjoint Sensitivity Analysis] \label{prop:adjoint_posterior}
	The posterior distribution $\beta | f_A \sim \bar{\mu}_A$ is given by 
	\begin{align*}	
		\bar{\mu}_A &= \mathcal{N}(\bar{a}_A, \bar{C}_A) \\
		\bar{a}_A &= a_A + C_A \adj{\mathcal{I}_A} (\mathcal{I}_AC_A\adj{\mathcal{I}_A})^{-1} (f_A - \mathcal{I}_Aa_A) \\
		\bar{C}_A &= C_A - C_A \adj{\mathcal{I}_A} (\mathcal{I}_AC_A\adj{\mathcal{I}_A})^{-1} \mathcal{I}_A C_A .
	\end{align*}
	The implied posterior distribution $\bar{\nu}_A$ is given by
	\begin{align*}
		\bar{\nu}_A &= \mathcal{N}(\bar{g}_A, \bar{G}_A) \\
		\bar{g}_A(p) &= -\mathcal{J}[p](\bar{a}_A) + \frac{\partial g}{\partial p}[U(p), p] \\
		\bar{G}_A(p,p') &= \mathcal{J}[p] C_A \adj{\mathcal{J}}[p']
	\end{align*}
\end{proposition}

Note that the form of the posterior over $\beta$ is essentially identical to the form of the posterior from \cref{prop:forward_posterior}, modulo the choice of information functionals and prior.
Next we will present some theoretical analysis of the forward and adjoint methods.

\subsection{Theoretical Analysis}

Our first theoretical result concerns a local error bound for the posterior mean in terms of the posterior covariance.
This result is a general result about conditional distributions of Gaussian process, and so is not specific to either the forward or adjoint method; as a result we adopt generic notation.

\begin{proposition}[Local error bound] \label{prop:local_error}
	Let $\mu = \mathcal{N}(a, C)$ be the prior, for $a \in \mathcal{H}_C$, and let $\bar{\mu} = \mathcal{N}(\bar{a}, \bar{C})$ be the posterior measure based on observations
	$\mathcal{I} u^\dagger = f$
	where $u^\dagger \in \mathcal{H}_C$, $\mathcal{I} \in (\mathcal{H}_C^\ast)^d$ and $f \in \reals^d$.
	Then we have that, for each $\mathcal{L} \in \mathcal{H}_C^*$
	\begin{equation*}
		\abs{\mathcal{L} \bar{a} - \mathcal{L}u^\dagger} \leq (\mathcal{L} \bar{C}\adj{\mathcal{L}})^\frac{1}{2} \norm{a - u^\dagger}_{C^{-1}} .
	\end{equation*}
\end{proposition}

The result from \cref{prop:local_error} is similar to results on error bounds in scattered data approximation with radial basis functions, such as in \citet{Wendland2004}.
The term $(\mathcal{L}\bar{C}\adj{\mathcal{L}})^\frac{1}{2}$ is analagous to the \emph{power function} \citet[Section 11.1]{Wendland2004}, but the focus in that work is on the case when both $\mathcal{L}$ and $\mathcal{I}_j$ are evaluation functionals.
In \citet[Chapter 16]{Wendland2004} each of these restrictions is relaxed, however the form of the power function derived in this case is more abstract than presented here.

Similar bounds appear in the literature on solution of PDEs by symmetric collocation with radial basis functions (see e.g. \citet[Section 16.3]{Wendland2004}, \citet{CockayneThesis, Cialenco2012}).
In these cases it is typically assumed that the $\tilde{\mathcal{I}_j}$ are evaluation functionals, so that the observations are point evaluations of the right-hand-side of the PDE, and that $\mathcal{L}$ is again an evaluation functional.
It is then possible to bound $(\mathcal{L} \bar{C} \adj{\mathcal{L}})^\frac{1}{2}$ in terms of the fill distance in the interior and on the boundary of the domain.
We have opted to make minimal assumptions on the form of the information operators and test functions in \cref{prop:local_error}, to avoid tying the result to a particular numerical method. 
Further note that the cited results only apply for fixed $p$ when performing sensitivity analysis for an elliptic PDE; as a global function of $(x, p)$ the sensitivity analysis equations may not be elliptic even when for fixed $p$ the underlying PDE is elliptic.

The next proposition provides theoretical guarantees for the setting when the solution $U(p)$ cannot be accessed directly, and instead a numerical estimate is provided by the map $\hat{U} : P \to \mathcal{U}$.
The natural way to provide such guarantees is by bounding the distance between the measure conditioned based on $U(p)$ to that based on $\hat{U}(p)$.
This is closely related to results that appear in \citet[e.g. Theorem 4.6]{Stuart:2010}, though the results presented therein assume that the two measures have a common dominating measure, which is not the case in the present setting.
A consequence of this is that the Hellinger metric, which is commonly used to measure distance in the space of probability measures in the field of uncertainty quantification, is not suitable here.

To proceed we introduce the $2$-Wasserstein metric, which is suitable for measures that are mutually singular.
Perhaps the most common way to define this metric is in terms of couplings of probability measures.
Let $\mu_1$ and $\mu_2$ be measures on some abstract normed space $\mathcal{V}$.
Let $\Gamma(\mu_1, \mu_2)$ be the set of couplings of $\mu_1$ and $\mu_2$, that is, the set of all Borel probability measures $\pi \in \mathcal{P}(\mathcal{V}\times \mathcal{V})$ with the property that $\pi(A \times \mathcal{V}) = \mu_1(A)$ and $\pi(\mathcal{V} \times A) = \mu_2(A)$ for each Borel set $A \subset \mathcal{V}$.
Then the $2$-Wasserstein metric \citep[Definition 6.1]{Villani2009} is given by
\begin{equation} \label{eq:wass_villani}
	W_2(\mu_1, \mu_2) = \left( \inf_{\pi \in \Gamma(\mu_1, \mu_2)} \int_{\mathcal{V}\times \mathcal{V}} \norm{v - v'}^2 \,\wrt\pi(v, v')
	\right)^\frac{1}{2} .
\end{equation}
We now proceed to state a generic result concerning robustness to approximation error, which will then be applied to the methods described in \cref{prop:forward_posterior} and \cref{prop:adjoint_posterior}.

\begin{proposition}[Robustness to Numerical Error] \label{prop:approx_error}
	Let $\mu = \mathcal{N}(a, C)$ be a Gaussian distribution with associated RKHS $\mathcal{H}$, for $a \in \mathcal{H}$ and $C : \mathcal{H} \to \mathcal{H}$ positive-definite.
	Assume that $\mathcal{I}$, $\hat{\mathcal{I}}$ are each bounded linear operators from $\mathcal{H}$ to $\reals^d$. 
	Let $\bar{\mu} = \mathcal{N}(\bar{a}, \bar{C})$ be the posterior measure based on observations
	$\mathcal{I} u^\dagger = f$
	where $u^\dagger \in \mathcal{H}$.
	Let $\hat{\mu} = \mathcal{N}(\hat{a}, \hat{C})$ be the same prior conditioned on observations
	$\hat{\mathcal{I}} u^\dagger = \hat{f}.$
	Then it holds that
	\begin{align*}
		W_2(\bar{\mu}, \hat{\mu}) &\leq 
			(C_{\mathcal{I},1} + C_{\mathcal{I},2}) \norm{\mathcal{I} - \hat{\mathcal{I}}}_{\mathcal{H}\to\reals^d} + C_f \norm{f - \hat{f}}_{\reals^d} + \mathcal{O}\left(\norm{\mathcal{I} - \hat{\mathcal{I}}}_{\mathcal{H}\to\reals^d}^2\right)
	\end{align*}
	where
	\begin{align*}
		C_{\mathcal{I},1} &= 
		\norm{C}_{\mathcal{H}} \left[
			\left(\norm{a}_\mathcal{H} + \norm{C}^{\frac{1}{2}}_{\mathsf{HS}}\right) 
			\left(
				\norm{G^{-1}\mathcal{I}}_{\mathcal{H}\to\reals^d} 
				+ \norm{\hat{G}^{-1} \hat{\mathcal{I}}}_{\mathcal{H}\to \reals^d} 
			\right) 
			+ \norm{G^{-1}f}_{\reals^d}
		\right]\\
		C_{\mathcal{I}, 2} &= 
		\alpha \norm{G^{-1}}_{\reals^d}^2 \norm{C}_{\mathcal{H}}\left(
			\norm{G^{-1}f}_{\reals^d} 
			+ \norm{\mathcal{I}}_{\mathcal{H}\to\reals^d}\norm{\hat{\mathcal{I}}}_{\mathcal{H}\to\reals^d} \left(\norm{a}_\mathcal{H} + \norm{C}^{\frac{1}{2}}_{\mathsf{HS}}\right)
		\right) \\
		C_f &= \norm{C}_{\mathcal{H}}\norm{G^{-1}\hat{\mathcal{I}}}_{\mathcal{H} \to \reals^d} \\
		\alpha &= \norm{\mathcal{I}C}_{\mathcal{H} \to \reals^d} + \norm{\hat{\mathcal{I}}C}_{\mathcal{H} \to \reals^d}
	\end{align*}
	and $G = \mathcal{I} C \adj{\mathcal{I}}$, $\hat{G} = \hat{\mathcal{I}} C \adj{\hat{\mathcal{I}}}$.
\end{proposition}

We next prove a corollary of this result which establishes a bound for the error in the posterior distribution for both forward and adjoint sensitivity analysis as a result of the need to use $\hat{U}(p)$ rather than having access to $U(p)$ directly.

\begin{corollary} \label{corr:robustness_specific}
	Assume that for each $p \in P$ there exists $\epsilon > 0$ such that
	\begin{align*}
		\Norm{\frac{\partial F}{\partial u}[U(p), p] - \frac{\partial F}{\partial u}[\hat{U}(p), p]}_{\mathcal{U} \to \mathcal{F}} &\leq \epsilon,\\
		\Norm{\frac{\partial F}{\partial p}[U(p), p] - \frac{\partial F}{\partial p}[\hat{U}(p), p]}_{P \to \mathcal{F}} &\leq \epsilon, \quad \text{and} \\
		\Norm{\frac{\partial g}{\partial u}[U(p), p] - \frac{\partial g}{\partial u}[\hat{U}(p), p]}_{\mathcal{U} \to \mathcal{G}} &\leq \epsilon.
	\end{align*}
	Further assume that the $\tilde{\mathcal{I}}_{F,j}$ and $\tilde{\mathcal{I}}_{A, j}$ are such that, for all $j = 1,\dots,d$
	\begin{align*}
		\norm{\tilde{\mathcal{I}}_{F,j}}_{\mathcal{F} \to \reals} < M < \infty, \\
		\norm{\tilde{\mathcal{I}}_{A,j}}_{\mathcal{F} \to \reals} < M < \infty .
	\end{align*}
	Lastly assume that $\norm{\quark}_{\reals^d} = \norm{\quark}_2$.

	Let $\hat{\mu}_F$ be the posterior distribution from \cref{prop:forward_posterior}, with $\hat{U}(\quark)$ substituted for $U(\quark)$.
	Likewise let $\hat{\mu}_A$ be the posterior from \cref{prop:adjoint_posterior} with the same substitution.
	Then we have
	\begin{align*}
		W_2(\bar{\mu}_F, \hat{\mu}_F) \leq (C_{\mathcal{I},1} + C_{\mathcal{I},2} + C_f)  M\epsilon\sqrt{d} + \mathcal{O}(\epsilon^2) \\
		W_2^2(\bar{\mu}_A, \hat{\mu}_A) \leq (C_{\mathcal{I},1} + C_{\mathcal{I},2} + C_f)  M\epsilon\sqrt{d} + \mathcal{O}(\epsilon^2)
	\end{align*}
\end{corollary}

\subsection{Comparison of Forward and Adjoint Approaches}

We conclude this section with a brief discussion of the relative merits of the forward and adjoint approaches, compared to the classical approach.

\paragraph{Choice of Method}
The forward approach requires the user to specify a prior on the parameter space; this is a space of dimension $\textup{dim}(P)$.
While the space in which the prior is placed for the adjoint problem is less directly connected to the derivative of interest, which might make eliciting a prior more challenging, in the finite-dimensional case reasoning about the correlation structure between the components of $\frac{\partial u}{\partial p}$ for the forward problem may also be challenging.
As a result, much as in classical sensitivity analysis, we are inclined to recommend the adjoint approach whenever $\textup{dim}(\mathcal{G}) < \textup{dim}(P)$, as will often be the case.
However if the user has strong prior information about the correlation structure between these components, the forward approach may still perform well.
Indeed, in the infinite-dimensional case such information is provided by knowledge about the smoothness of the function $p$.

\paragraph{Experimental Design}
\cref{prop:forward_posterior,prop:adjoint_posterior} each allow the user to construct a global model for the required derivatives.
However in order to perform inference globally, one requires a set of points in $P$ with which to construct the posterior.
Both the forward and the adjoint approach suffer from the curse of dimensionality in this respect, since Gaussian processes typically require such designs to be ``space-filling''\footnote{
	Since, typically, the rate of convergence of Gaussian processes with this type of information depends on the ``fill distance'', i.e.~the maximum distance of any point in the space to a design point. 
	See e.g. \citet{Wendland2004}, or \citet{Cialenco2012} in the context of PDEs.
}, and if $P$ is high-dimensional constructing a space filling design will be equally prohibitive in either mode.
However in the present paper we focus on application of these methods within iterative optimisation algorithms, so that rather than requiring a space-filling design we only require good estimates of the gradient along the path in parameter space followed by the optimiser.
This will be discussed in detail in the next section.

\section{Optimisation and Probabilistic Sensitivity Analysis} \label{sec:optimisation}

We now explore a potential application of probabilistic local sensitivity analysis, as a way to provide approximations of gradients in optimisation algorithms.
As a starting point we will consider the most fundamental of gradient-based optimisation algorithms, gradient descent \citep{Curry1944}.
In \cref{sec:grad_descent} we briefly recall the GD algorithm.
In \cref{sec:probabilistic_grad_descent} we describe how probabilistic gradients can be incorporated into the algorithm.
Then, in \cref{sec:applications} we explore the use of this approach in two applications.

\subsection{Gradient Descent} \label{sec:grad_descent}

We now describe the GD algorithm.
GD is in many respects a prototypical gradient-based optimisation method, making it a natural starting point for studying the integration of probabilistic gradients into such algorithms.
In GD the goal is to compute a (local) minimiser $p^*$ of a function $g(p)$.
To accomplish this a sequence of points $(p^n)$, $p^n \in P, n\in \mathbb{N}$ is generated iteratively starting from some user-defined initial point $p^0$ and advancing according to 
\begin{equation*}
	p^{n+1} = p^n - \gamma^n \frac{\wrt g}{\wrt p}(p^n)
\end{equation*}
where $\gamma^n$ is a parameter of the method known as the \emph{step size} or \emph{learning rate}.
Under specific conditions on $f$ and $\gamma^n$ it can be shown that $p^n \to p^*$ (again, a local minimiser) as $n\to\infty$; see \citet[Section 3.2]{Nocedal2006} for further details.
GD is presented as an algorithm in \cref{alg:grad_descent} in the supplement.

There are various methods for choosing the parameter $\gamma_n$.
Since the focus of this work is on the performance when $\frac{\wrt g}{\wrt p}$ is replaced by the probabilistic gradients introduced in \cref{sec:probabilistic}, we will use a probabilistic version of the backtracking line search method described in \citet[Algorithm 3.1]{Nocedal2006}, based on the method described in \citet{Mahsereci2015}.

\subsection{Gradient Descent with Probabilistic Gradients} \label{sec:probabilistic_grad_descent}

We now discuss a probabilistic modification of GD.
Heuristically the approach followed is to replace the computation of $\frac{\wrt g}{\wrt p}$ with a probabilistic gradient obtained from either \cref{prop:forward_posterior} or \cref{prop:adjoint_posterior}; to simplify the exposition we will describe the former, but the approach is essentially identical in the latter.
The approach is presented as an algorithm in \cref{alg:p_grad_descent}.
Essentially, we begin with a prior $\mu_F$ which is projected to $\nu_F$ as described in \cref{prop:forward_posterior}.
We then construct a sequence of random variables $(X_F^n)$, where $X_F^0$ has law $\nu_F$, by sequentially updating this prior with information collected over the course of the optimisation.
This provides a posterior distribution over the gradient which is used in place of $\frac{\partial g}{\partial p}$ in GD.
The principal advantages, illustrated in \cref{sec:applications}, are that (i) for each value of $p^n$, one can often obtain an approximation of $\frac{\wrt g}{\wrt p}$ that is sufficiently accurate for the purposes of taking a gradient step, at a lower cost than that of computing $\frac{\wrt g}{\wrt p}$ directly, and (ii) since the posterior is defined over the entire parameter space, for some values of $p^n$ \emph{no} inversion problem must be solved to advance the gradient descent.

There are two main issues to address.
The first is that that it is well-established in the literature on stochastic gradient descent that line-search algorithms such as the \textsc{bls} routine are not robust to inaccurate gradients.
This is discussed in \citet{Mahsereci2015}.
Since the gradients we propose to use in this work are also inaccurate, an alternative line-search strategy for selecting the step sizes $\gamma^n$ must be adopted in the probabilistic case.
Borrowing from the literature on stochastic gradient descent, our proposed approach incorporates ideas from the probabilistic line search of \citet{Mahsereci2015} into the backtracking line search from \citet[Section 3.2]{Nocedal2006}.
The \textsc{pls} routine is described in \cref{alg:pls}.

A second issue is that if $X_F^n$ is not sufficiently accurate, the step size $\gamma^n$ found by the probabilistic line search will be selected to be below the tolerance $\epsilon$, causing the algorithm to terminate.
To address this we propose to couple the computation of $\gamma^n$ with the calculation of the gradient, as described in \textsc{probjac} within \cref{alg:p_grad_descent}. 
Once the tolerance has been achieved, we calculate the step size $\gamma$ according to a probabilistic version of backtracking line search that will be described presently.
If $\gamma$ is above the tolerance the procedure returns the current gradient estimate, along with the posterior distribution and the step size; otherwise, the tolerance $\delta$ is reduced and the conditioning procedure is repeated.
This continues until delta is below some minimum value $\delta_{\min}$, at which point convergence is accepted.

\begin{algorithm}
\caption{
	Probabilistic version of gradient descent.
	The routines \textsc{metric} and \textsc{info} are problem specific and must be supplied by the user, with the former assessing the distribution of the currently computed posterior distribution to determine whether it is sufficiently narrow to accept it as a valid gradient and the latter supplying information, iteratively, based on the current distribution and location.
	The routine \textsc{condition} implements \cref{prop:forward_posterior}.
	\textsc{pls} is the probabilistic version of the Armijo line search, and is given in \cref{alg:pls}.
	Of the new parameters, $\delta$ reflects how much accuracy is demanded of the posterior at each iteration, $\delta_{\min}$ specifies a maximum level of accuracy to protect against numerical instabilities resulting from large Gram matrices in \textsc{condition}, and $\tau$ describes how rapidly $\delta$ is reduced when a valid descent direction cannot be found.
} \label{alg:p_grad_descent}
\begin{algorithmic}[1]
	\Procedure{\textsc{pgd}}{$p^0, g, \mu^0_F, \epsilon, \delta, \delta_{\min}, \tau_1$}
		\State Compute $\nu_F^0$ from $\mu_F^0$ and let $X_F^0$ be the random variable with law $\nu_F^0$
		\For{$n = 1,2,\dots$}
			\State $s_F^n, X_F^n, \gamma_n \gets \textsc{probjac}(X_F^{n-1}, g, p^{n-1}, \epsilon, \delta, \delta_{\min})$
			\If{$\gamma^n < \epsilon$}
				\State \Return $p^{n-1}$
			\EndIf
			\State $p^n \gets p^{n-1} + \gamma^n s_F^n$
		\EndFor
	\EndProcedure

	\Procedure{\textsc{probjac}}{$X, g, p, \epsilon, \delta, \delta_{\min}, \tau_1$}
		\While{$\delta > \delta_{\min}$}
			\While{$\textsc{metric}(X) > \delta$}
				\State $\mathcal{I}, f \gets \textsc{info}(X, p)$
				\State $X \gets \textsc{condition}(X, \mathcal{I}, f)$
				\State $s \gets -\mathbb{E}(X(p)) / \norm{\mathbb{E}(X(p))}_2$
    			\State $\gamma \gets \textsc{pls}(p, g, X)$
    			\If{$\gamma < \epsilon$}
    				\State $\delta \gets \tau_1\delta $
    			\Else
    		        \State \Return $s, X, \gamma$
    			\EndIf
			\EndWhile
		\EndWhile
	\EndProcedure
\end{algorithmic}
\end{algorithm}

\begin{algorithm}
\caption{Probabilistic line search algorithm. 
This is essentially a modification of the backtracking line search described in \citet[Algorithm 3.1]{Nocedal2006} to account for the fact that the gradient is a random variable rather than a constant.
The parameters $p, g$ and $X$ are the parameter value, objective function and current posterior, respectively. 
The remaining parameters control the behaviour of the algorithm; we have specified sensible defaults for these and assume those defaults are used throughout the text.
$\tau_2$ controls how rapidly $\gamma$ is decreased, while $c$ controls how large a reduction in the objective function is required when a step is taken in the chosen direction and $P^\textup{crit}$ is the probability with which this reduction must be achieved. $\gamma$ and $\gamma^{\textup{min}}$ control the initial and minimum values of $\gamma$ respectively.
}\label{alg:pls}
\begin{algorithmic}
	\Procedure{\textsc{pls}}{$p, g, X; \;\tau_2=0.5, c=0.5, P^\textup{crit}, \gamma=1, \gamma^{\textup{min}}=10^{-10}$}
		\State $s \gets - \mathbb{E}(X) / \norm{\mathbb{E}(X)}_2$
		\While{$\gamma > \gamma^{\textup{min}}$}
			\State $p_\gamma \gets p + \gamma s$
			\If{$g(p_\gamma) > g(p)$}
				\State \textbf{continue}
			\EndIf
			\State $Z \gets -c\gamma X^\top s$
			\If{$\mathbb{P}(Z > g(p_\gamma) - g(p)) < P^\textup{crit}$}
				\State \Return $\gamma$
			\EndIf
			\State $\gamma \gets \tau_2\gamma$
		\EndWhile
	\EndProcedure
\end{algorithmic}
\end{algorithm}

\subsubsection{Discussion} \label{sec:algorithm_discussion}
We now provide some important remarks about the algorithm presented above.

\paragraph{Choice of Direction}

The direction chosen in \cref{alg:p_grad_descent} at each iteration is the posterior mean.
A natural alternative would be to instead \emph{sample} a direction from the posterior distribution.
This requires only minor modification of the above algorithm, but empirically was found to perform slightly worse in general; consequently we have opted to use the posterior mean as the descent direction.

\paragraph{Recycling Information}
Note that the gradient here is computed based on information collected at all points $p_F^1,\dots,p_F^n$, i.e.~based on a \emph{global} model for the gradient as a function of $p$.
Since the sequence $(p_F^n)$ will increasingly concentrate in a region of $p_*$ as $n$ increases, one expects that the prior $\bar{\mu}_F^{n-1}$ will be an increasingly accurate predictor for the gradient $\frac{\wrt g}{\wrt p}(p_F^n)$ as $n$ increases.
This means that once some computational effort has been expended to obtain a relatively accurate gradient, it is possible for \textsc{probjac} to perform many further iterations based on this gradient without needing calls to \textsc{condition}, as we shall see in \cref{sec:applications}.

\paragraph{Linearly Independent Information}
A global model introduces some additional burden to ensure that $\mathcal{I}_F^n$ is linearly independent of $\mathcal{I}_F^1,\dots,\mathcal{I}_F^{n-1}$, both to maximise the amount of new information obtained at each $p^n$ and to ensure that the linear system that must be solved to compute the posterior does not become singular.
Thus, $\textsc{info}$ must be carefully designed to ensure that the information returned is not too highly correlated with information already observed.

\paragraph{Computational Cost}
To compute the posterior distributions from \cref{prop:forward_posterior} and \cref{prop:adjoint_posterior}, it is necessary to compute the matrix $M = (\mathcal{I} C \adj{\mathcal{I}})^{-1} \mathcal{I} C$ by solving the linear system
$\mathcal{I} C \adj{\mathcal{I}} M = \mathcal{I} C$.
To accomplish this one typically computes a Cholesky factorisation of $\mathcal{I} C \adj{\mathcal{I}}$, which becomes computationally intensive once many information functionals have been collected.
However, we note that the sequential nature of the algorithm proposed is such that, rather than recomputing the full factorisation at each iteration of \textsc{probjac}, one can use an updating formula for the factorisation such as presented in \citet[Appendix B]{Osborne:2010}; this is described in detail in \cref{sec:posterior_updates}.
In brief, one must only compute the Cholesky factorisation of a smaller matrix, whose dimension is only the same size as the dimension of the \emph{new} information, which naturally dramatically reduces the cost of computing the probabilistic gradients.

The other factor that influences the cost is the size of $\mathcal{I}C$, and since this defines \emph{how many} linear systems must be solved, it may be that ultimately the cost of assembling the posterior $\bar{\mu}_F^n$ exceeds than that of simply computing $\frac{\wrt g}{\wrt p}(p_F^n)$ despite the efficient updating formula for the factorisation.
Thus in practise we propose that the \textsc{probjac} is used only to perform the initial iterations, and that when the method is determined to be close to the truth, or the cost of constructing the posterior is too great, we revert to classical GD to complete the optimisation.
In \cref{sec:applications} we adopt the crude rule of thumb that \textsc{probjac} is terminated when the dimension of $\mathcal{I} C \adj{\mathcal{I}}$ exceeds $10,000$, though this is never exceeded in practise for one of the two examples examined.
In future work more sophisticated switching schemes will be explored.

\paragraph{Choice of Metric}

The routine \textsc{metric} must assess whether the posterior distribution at a particular iteration is sufficiently accurate for the probabilistic gradient to be accepted as a valid direction for the gradient descent.
To determine this we focus on the width of the posterior covariance, and in this work we exclusively use the square-root of the trace of the posterior covariance, $\sqrt{\textup{trace}(\bar{G}_F)}$ as a proxy for the width.
An exploration of other choices is not expected to affect the performance of the algorithm dramatically, and is reserved for future work.

\paragraph{Choice of Information Functionals}
Lastly, we note that we have not yet discussed the selection of information functionals in $\textsc{info}$.
We expect this to be highly problem dependent.
We make a proposal in the next section that appears to be well adapted to the two examples presented therein, but do not expect that there exists a unique optimal choice of information for all settings.

\section{Applications} \label{sec:applications}

In this section we apply \cref{alg:p_grad_descent} to compute the maximum \emph{a-posteriori} (MAP) point in Bayesian inversion problems for two problems.
In \cref{sec:fitzhugh_nagumo} we seek to infer a small number of parameters of an ODE using the forward approach, and in \cref{sec:gwf} inference of a larger number of parameters of a challenging PDE using the adjoint approach.

\subsection{FitzHugh---Nagumo Model} \label{sec:fitzhugh_nagumo}

As a first example we examine the problem of inferring the parameters for the Fitzhugh---Nagumo model \citep{FitzHugh1961}, a nonlinear oscillatory ODE.
Since this problem has four parameters, we use the forward approach from \cref{sec:probabilistic:forward}.

\subsubsection{Problem Definition}

The equations that define the FitzHugh---Nagumo model are
\begin{align*}
	\frac{\wrt v}{\wrt t} &= v - \frac{v^3}{3} - w + I &\qquad
	\frac{\wrt w}{\wrt t} &= \frac{v + a - bw}{\tau}
\end{align*}
where $a, b, I, \tau \in \reals_+$ are parameters of the model.
We concatenate the parameters as $p = [I, a, b, \tau]^\top \in \reals^4 = P$.
The solution to this system of ODEs for $p^*=[0.5, 0.8, 0.7, 12.5]^\top$ is shown in the supplement in \cref{fig:fhn:solution}, while sensitivities are displayed in \cref{fig:fhn:dv_dp,fig:fhn:dw_dp}.
The solution space $\mathcal{U}$ is a space of once differentiable functions $u : D \to \reals^2$, where $D = [0, T]$ for some $T > 0$.
A reformulation of this problem in terms of the constraint function $F(u,p)$ can be found in \cref{sec:fitzhugh_nagumo:supplement}, along with the form of its derivatives $\frac{\partial F}{\partial u}$, $\frac{\partial F}{\partial p}$.

To set up the inference problem we generated data for true parameter values $p^*$ by evaluating the $v(t_i^\textup{data}; p^*)$ at times $t_i^\textup{data} = i$, $i = 1,\dots,20$.
These locations are distinguished as dashed gray lines in \cref{fig:fhn} in the supplement.
Observations were then corrupted with centred Gaussian noise with standard deviation $10^{-2}$, i.e.
$y_i = v(t_i^\textup{data}; p^*) + \xi_i$
where $\xi_i \sim \mathcal{N}(0, 10^{-2})$ IID.
The prior distribution over the parameters was set to be log-Gaussian with mean $m_p=[1, 1, 1, 10]^\top$ and covariance $I$.
The objective function is twice the negative logarithm of the likelihood multiplied by the prior, and is thus given by
\begin{equation*}
	g(p) := \frac{1}{\gamma^2}\sum_{i=1}^M (v(t_i^\textrm{data}; p) - y_i)^2 + (\log(p) - \mu)^\top \Sigma^{-1} (\log(p) - \mu)
\end{equation*}

\subsubsection{Probabilistic Gradient Descent}

To apply the probabilistic gradient descent algorithm from \cref{alg:p_grad_descent} we must first specify the prior over $\frac{\wrt U}{\wrt p}$.
Since the parameter space is four-dimensional and $\mathcal{U}$ is a space of vector-valued functions, formally $\frac{\wrt U}{\wrt p}$ is $\reals^{2 \times 4}$-valued.
For convenience, we place a prior on $X : D \times P \to \reals^8$, and form $\frac{\wrt U}{\wrt p}$ as
\begin{equation*}
	\frac{\wrt U}{\wrt p} = \begin{bmatrix}
		X_{1:4}^\top \\
		X_{5:8}^\top
	\end{bmatrix}
\end{equation*}
where $X_{i:j}$ denotes components $i$ to $j$ of $X$.
Noting that the posterior covariance is independent of the data, we assume an independent and identical prior over each column of $\frac{\wrt U}{\wrt p}$, so that the inference is identical but for the distinct right-hand-side for each component of $p$ in the posterior mean of \cref{prop:forward_posterior}.

Since the initial condition is independent of $p$, this prior was taken to be $X \sim \mathcal{N}(\bm{0}, k)$ where
\begin{align*}
	k((t, p), (t', p')) &= 
		C 
		 q(t) q(t') k_{5/2}([t, p]^\top, [t', p']^\top \sigma, L) \\
	k_{5/2}(r, r'; \sigma, L) 
	&= \sigma^2 \left( 
		1 + \sqrt{5} d(r, r'; L) + \frac{5}{3} d(r, r'; L)^2
	\right)
	\exp\left(
		-\sqrt{5}d(r, r'; L)
	\right) \numberthis\label{eq:matern_52} \\
	d(r, r'; L) &= \sqrt{r^\top L^{-1} r'} \\
	C &= \begin{bmatrix}
		1 & \rho \\
		\rho & 1
	\end{bmatrix} \\
	q(t) &= t .
\end{align*}
Multiplication by the linear functions $q(t)$ ensures that there is no uncertainty at $t=0$, where the sensitivity is known to be zero.

The kernel $k_{5/2}$ in \cref{eq:matern_52} is a member of the Matérn family \citep[Section 4.2]{Rasmussen2005} and is the covariance kernel for a prior over functions with at least two continuous derivatives.
To ease computation the length-scale matrix $L$ was selected to be diagonal, $L = \textup{diag}(\ell)$ for $\ell \in \reals^6$.
This parameter was further restricted to $\ell = [\ell_x \bm{1}_2, \ell_p \bm{1}_4]$ where $\ell_x, \ell_p \in \reals$.
The scalars $\sigma, \ell_x$ and $\ell_p$ were then selected by maximising the marginal likelihood of an initial candidate design \citep[see e.g.][Section 5.4]{Rasmussen2005}.
This was obtained by sampling a set of candidate parameters $p_i^\textup{calib}$, $i=1,\dots,5$ from the prior over the parameters and defining the corresponding evaluation functionals $\tilde{I}_{ij} = \delta[i]$, $i=1,\dots,20$ (i.e.~using equally spaced points inside the spatial domain).
The parameter $\rho$, which describes the degree of prior covariance between the components $u_1 = v$ and $u_2 = w$, was fixed to 0.5.

For this problem it was convenient to restrict the information functionals to be evaluation functionals, i.e.~$\mathcal{\tilde{I}}_i = \delta[t_i^\textup{info}]$.
The points $t_i^\textup{info}$ were restricted to a fine grid of $1000$ points in $(0,T]$, denoted $t_1^\textup{info},\dots,t_{1000}^\textup{info}$.
To choose the next information functionals at iteration $n$ within the function $\textsc{info}$ in \cref{alg:p_grad_descent}, we choose new conditioning locations within this set by attempting to minimise a heuristic based on the fill distance which often appears as an upper bound in Gaussian process regression problems.
To be specific, we begin by constructing an augmented point set:
\begin{equation*}
	z_{ij} = \begin{bmatrix}
		t_i^\textup{info} \\ p^j
	\end{bmatrix} .
\end{equation*}
for $j=1,\dots,n$ denoting the iteration number in \textsc{probjac} and $p^j$ the corresponding parameter value for that iteration.
The information functionals were then selected to be the $\tilde{\mathcal{I}}_j$ for which the distance between $z_{in}$ and $z_{i'j}$, $i,i'=1,\dots,1000$, $j=1,\dots,n-1$, is maximised.

\subsubsection{Results}

\begin{figure}
	\centering
	\begin{subfigure}{0.49\textwidth}
		\includegraphics[width=\textwidth]{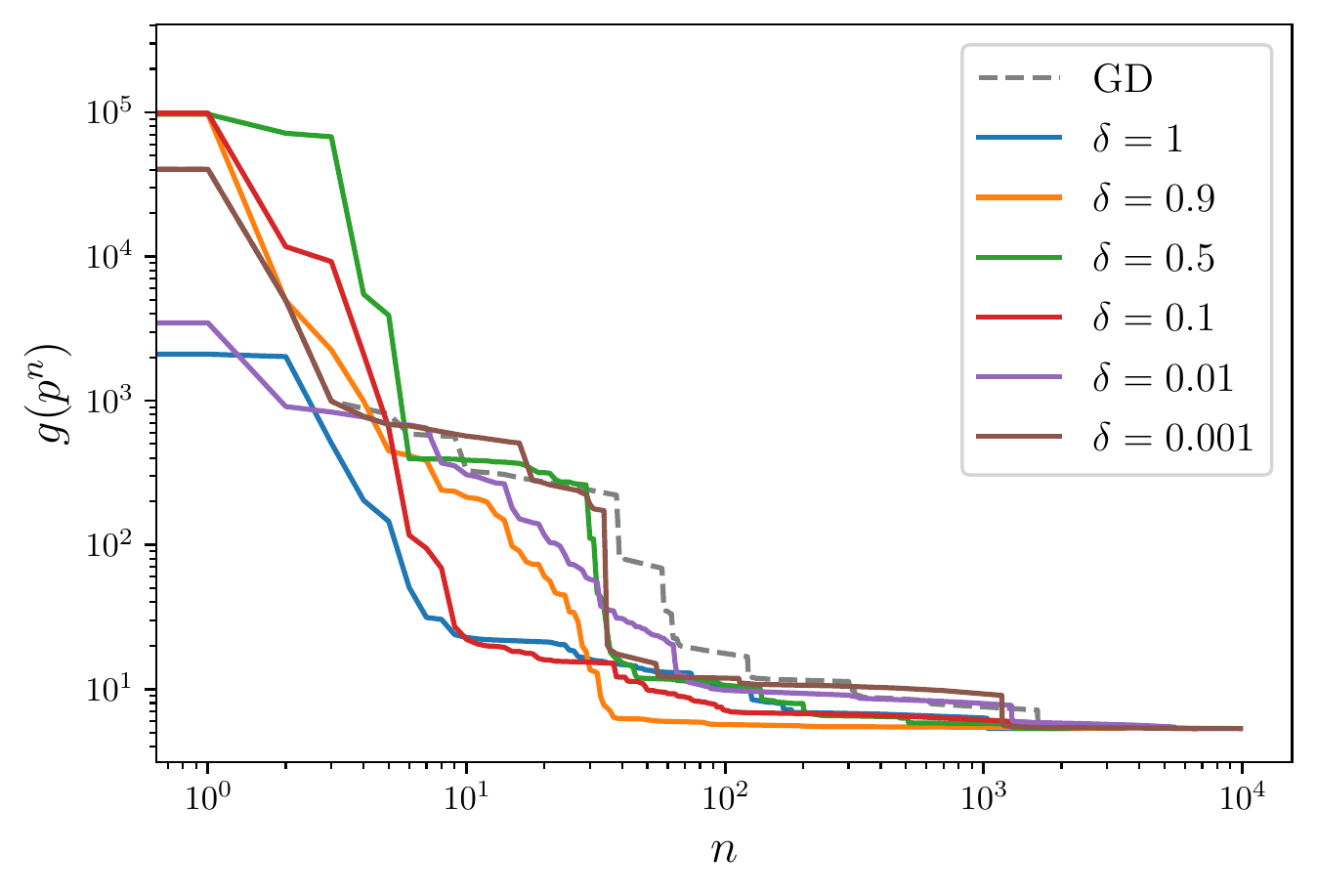}
		\caption{$g(p^n)$} \label{fig:fn_results_g}
	\end{subfigure}
	\begin{subfigure}{0.49\textwidth}
		\includegraphics[width=\textwidth]{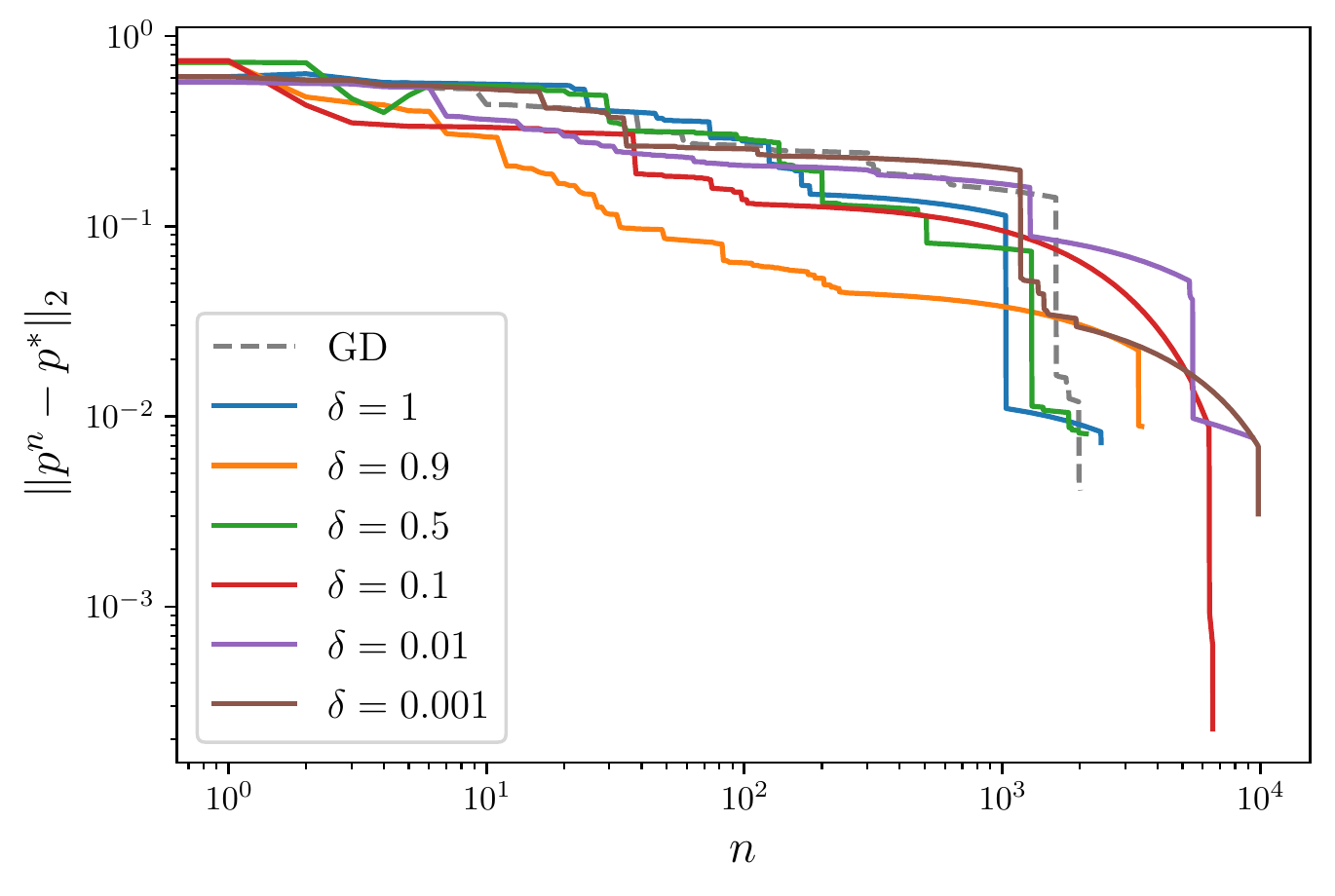}
		\caption{$\norm{p^n - p^*}_2$}  \label{fig:fn_results_distance}
	\end{subfigure}

	\begin{subfigure}{0.49\textwidth}
		\includegraphics[width=\textwidth]{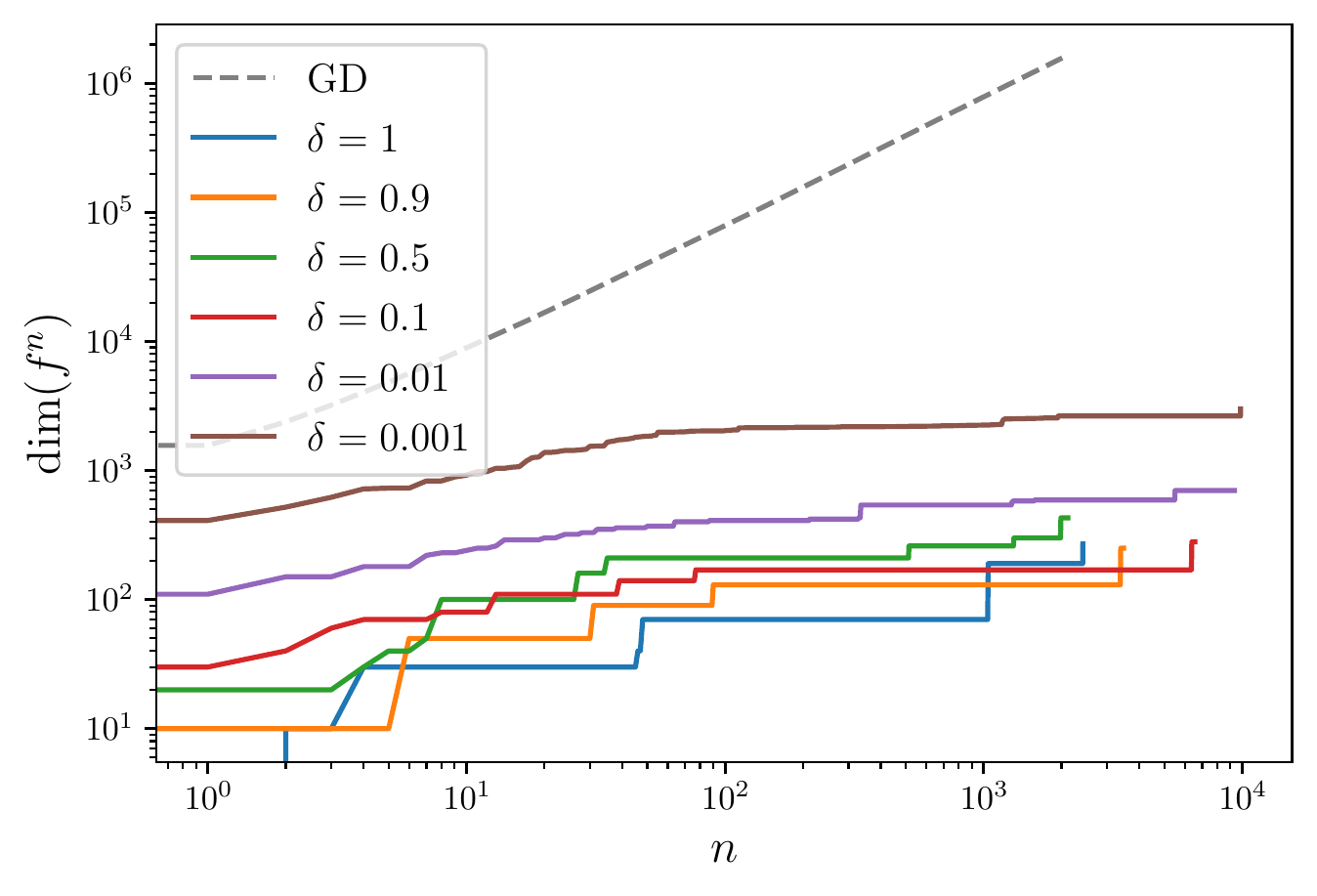}
		\caption{$\textup{dim}(f^n)$}  \label{fig:fn_results_dim}
	\end{subfigure}
	\caption{Performance of the probabilistic gradient descent algorithm on the FitzHugh---Nagumo model described in \cref{sec:fitzhugh_nagumo} as the parameter $\delta$, which roughly controlling the accuracy demanded of the probabilistic gradient estimate, is varied.
	Here $n$ is the iteration number.
	\cref{fig:fn_results_dim} shows the value of the quantity of interest $g$, in this case the value of the negative log-target in a Bayesian inference problem described in \cref{sec:fitzhugh_nagumo}.
	\cref{fig:fn_results_distance} shows the distance from the parameter at iteration $n$ to the true MAP point.
	\cref{fig:fn_results_dim} shows the dimension of the matrix inversion problem that was solved in order to compute the posterior distribution.
	}
	\label{fig:fn_results}
\end{figure}

The paths taken by the probabilistic optimiser are contrasted with classical gradient descent in \cref{fig:fn_results}.
\cref{fig:fn_results_g} shows the value of $g(p^n)$, while \cref{fig:fn_results_distance} shows the distance from the minimum obtained by gradient descent.
All of the methods were started from the initial parameter value $p_0 = m_p$, and the GD tolerance was set to $\epsilon = 10^{-6}$.
The threshold $\delta$ was varied from $1$ (representing a high level of allowed error in the posterior gradient estimate) to $0.001$ (representing a low level of allowed error).
In each case $\delta_{\min}$ was set to $10^{-6}$.
For $\delta=0.9, 0.5, 0.1$ the performance of the probabilistic approach is initially worse, as expected, though as the iterates near $p^*$ the performance of the probabilistic approaches improves.
Interestingly, for $\delta=1$ and $\delta=0.01$ the probabilistic approach actually seems to initially converge \emph{faster} than the classical approach.
This should not generally be expected, though we note that since the GD directions have no particular optimality properties nothing prevents an approximate method from achieving faster convergence.

\cref{fig:fn_results_dim} tracks the amount of data collected (i.e.~the size of $f_F^n$) as a function of the iteration number.
This exhibits the expected behaviour of increasing inversely proportional to $\delta$.
However it is noteworthy that even in the strictest case, $\delta=0.001$, only $3000$ evaluations of $\frac{\partial F}{\partial p}$ are required over the course of $9840$ iterations to perform almost as well as as gradient descent.
For context computing the gradient $\frac{\wrt g}{\wrt p}$ using the \texttt{DOP853} algorithm \citep[Section II]{Hairer1993} method as implemented in \texttt{scipy} required an average of $781$ evaluations of $\frac{\partial F}{\partial p}$ \emph{per iteration} of gradient descent, with a total of over $1.5$ million evaluations over the course of the $2013$ iterations performed with exact gradients.
While $\frac{\partial F}{\partial p}$ is cheap to evaluate in this example, in a setting in which this was a bottleneck it is clear that the probabilistic method would be preferable.
Further note that while $3000$ evaluations of $\frac{\partial F}{\partial p}$ were required, as noted in \cref{sec:algorithm_discussion} this does not translate directly to inversion of a $3000\times 3000$ Gram matrix, as the updating formula for Cholesky factorisations was exploited.

\subsection{Groundwater Flow Model} \label{sec:gwf}

We now consider a linear PDE that describes the steady-state flow of fluid through a porous medium.
In this section the parameter is formally function-valued.
Since after discretisation its dimension can be large, the adjoint approach is adopted.

\subsubsection{Problem Definition}
For a fixed value of the parameter $p$, the forward model is given by
\begin{align*}
	-\nabla \cdot \left( p(x) \nabla u(x)\right) &= 0 & x &\in D \\
	u(x) &= x & x_2 &= 0 \\
	u(x) &= 1-x & x_2 &= 1 \\
	\frac{\partial u}{\partial x_1} &= 0 & x_1 &= 0 \textrm{ or } 1
\end{align*}
Here the domain $D = [0,1]^2$ and $p : D \to \reals$.
We assume that $p(x) > 0$ for all $x \in D$.

The solution $u(x)$ was obtained by discretising the domain above with FEM on a fine triangular meshing of the unit square based on a grid of $32\times 32$ points using piecewise-linear basis functions.
The mesh is depicted in \cref{fig:gwf_solution:mesh} in the supplement, and the discretisation results in a finite-dimensional approximation of the solution $u(x)$ in with $1089$ degrees of freedom.
The solution to the PDE above for the parameter value $p = 1$ is depicted in \cref{fig:gwf_solution:soln}, again found in the supplement.

The parameter is defined to be piecewise constant over supersets of the cells of this mesh, defined by grouping the cells based on a subdivision of the domain into squares.
For a parameter $N > 1$ these are obtained by placing down a regular grid of $N^2$ points, with $N+1$ equispaced points along each axis.
The points of this grid form the vertices of the $N^2$ parameter cells.
In \cref{fig:gwf_solution:mesh}, the parameter cells for $N=4$ are surrounded by green lines.

To construct the inverse problem, we use a Gaussian prior $p \sim \mathcal{N}(\mu, \Sigma) =: \pi(p)$ with $\mu = 5\bm{1}_{N^2}$, where $\bm{1}_{N^2}$ here denotes the vector of ones in $\reals^{N^2}$.
Letting $x_i^\textup{param}$ denote the centroid of cell $i$ according to some arbitrary ordering of the cells, $i=1,\dots,N$, the covariance is given by $\Sigma_{ij} = k(x_i^\textup{param}, x_j^\textup{param})$, where $k$ is the Matérn $5/2$ kernel given in \cref{eq:matern_52}, with amplitude and length-scale each set to 1.
The data-generating parameter $p^*$ was sampled from the prior over $p$.
To define the likelihood, data was obtained by taking direct measurements of the solution $u(x; p^*)$ at locations $x_1^\textrm{data},\dots x_M^\textrm{data}$ where $M=25$ the $x_j^\textrm{data}$ are the nearest mesh points to points on a regular $5 \times 5$ grid starting at $(0.1, 0.1)$ and ending at $(0.9, 0.9)$.
The points of this grid are shown in \cref{fig:gwf_solution:soln} in the supplement.
Let $\tilde{y} \in \reals^M$ be the vector with $\tilde{y}_j = u(x_j^\textrm{data}; p^*)$.
These points were corrupted with IID Gaussian noise $\xi_j \sim \mathcal{N}(0, \gamma)$, $\gamma=0.01$ $j=1,\dots,M$ to obtain data $y = \tilde{y} + \xi$.
Denoting the likelihood by $\pi(p | y, u)$ with dependence on $u$ emphasised, the QoI for gradient descent was then given by $g(u, p) = -2\log \pi(p | y, u) \pi(p)$, i.e.
\begin{equation} \label{eq:g_gwf}
	g(p) := \frac{1}{\gamma^2}\sum_{i=1}^M (u(x_i^\textrm{data}; p) - y_i)^2 + (p - \mu)^\top \Sigma^{-1} (p - \mu)
\end{equation}

\subsubsection{Probabilistic Gradient Descent}

To test the algorithm described in \cref{sec:probabilistic_grad_descent} we attempt to compute the MAP point of the posterior distribution for the inverse problem described above.
Owing to the potentially high dimension of the problem to be solved, the adjoint approach was used.
For the prior we used $\beta \sim \mathcal{N}(0, k)$, where $k$ is given by
\begin{equation*}
	k((x, p), (x', p')) = k_{52}((x, p), (x', p'); \sigma, \ell) q(x_2) q(x_2').
\end{equation*}
Here $q(x) = 1 - (2x-1)^2$, so that $q(0) = q(1) = 0$, ensuring that the relevant boundary condition is encoded in the prior since we note that the boundary conditions do not depend upon $p$.
Thus, the prior is formally over functions from $\reals^{N^2 + 2}$ to $\reals$, though since the problem has been discretised with the finite-element method the discretised prior is finite-dimensional.
Strictly speaking to project the prior into the finite-element space requires computing integrals of the form $\int k(x, x') \phi_j(x) \wrt x$ for $j=1,\dots,1089$, however since these integrals do not generally have a closed-form we opt to approximate them as 
$\int k(x, x') \phi_j(x) \wrt x \approx k(x_j, x')$
where $x_j$ is the nodal point corresponding to the basis function $\phi_j$.

For the parameters of the prior, a separate constant length-scale was assigned to the spatial variables and the parameters, denoted $\ell_x$ and $\ell_p$ respectively, i.e.~$\ell = [\ell_x \bm{1}_2, \ell_p\bm{1}_{N^2}]$.
The amplitude $\sigma$ and the length-scale $\ell_p$ were again selected by maximising the marginal likelihood of these parameters given a candidate design obtained again by sampling a set of candidate parameters $p_i^\textup{calib}$, $i=1,\dots,10$, from the prior over parameters, and choosing corresponding information functionals $\tilde{\mathcal{I}}_{ij}u = \int_D u(x) \phi_j(x) \wrt x$. 
Here the $\phi_j$ are the finite element basis functions corresponding to the nearest mesh points to a regular $10\times 10$ grid of points within $D$, with basis functions on the top and bottom boundaries excluded.

For the remaining parameter, $\ell_x$, we note that since in \cref{eq:g_gwf} $g$ depends only on the value of $u$ at the points $x_i^\textrm{data}$, we therefore have that $\frac{\partial g}{\partial u}$ is zero everywhere but at these locations.
Since this function is so rough, it is impossible to infer the spatial length-scale $\ell_x$ from evaluations of it.
As a result, we opted to fix $\ell_x = 0.2$, based on the observed smoothness of the solution to the adjoint equations.

For the information functionals we selected $\tilde{\mathcal{I}}_j u = \int u(x) \phi_j(x) \wrt x$, i.e.~projection against the $j$\textsuperscript{th} finite element basis function.
This is straightforward to implement since after discretisation it is simply projection against the canonical basis vector $\bm{e}_j^\top$.
The function $\textsc{info}$ was implemented similarly to in \cref{sec:fitzhugh_nagumo}, with the fine grid of points now consisting of the mesh locations which the basis functions correspond to, again excluding points on the top and bottom boundaries.
However, to ensure that the information $f$ is nonzero, we enforce that when $\textsc{metric}(X^n) > \delta$, the first locations to be conditioned upon are those basis functions corresponding to $x_i^\textrm{data}$.

\subsubsection{Results}

\begin{figure}
    \centering
	\begin{subfigure}{0.7\textwidth}
		\includegraphics[width=\textwidth]{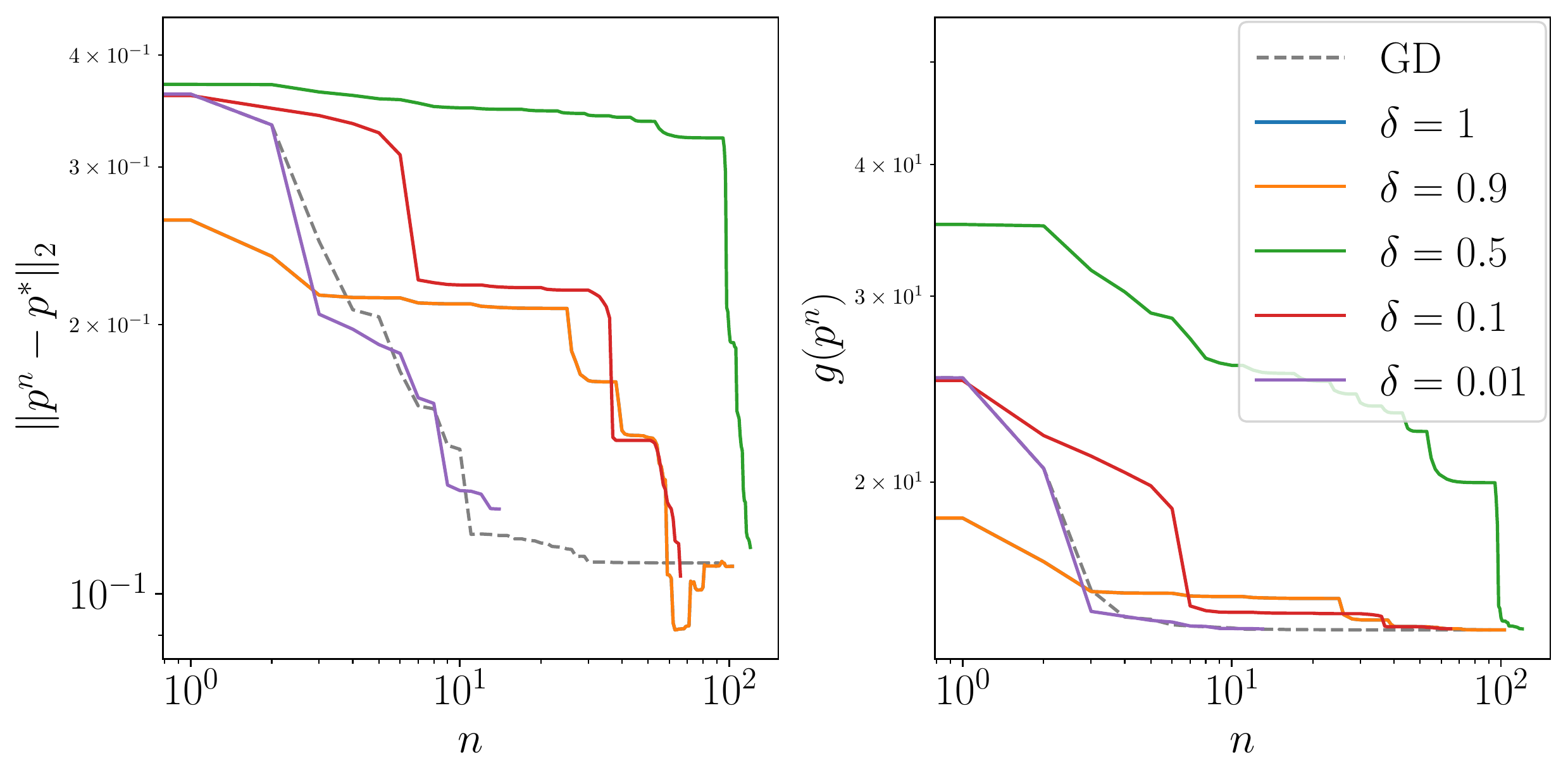}
		\caption{$N=2$; $\textup{dim}(P)=4$}
	\end{subfigure}
	
	\begin{subfigure}{0.7\textwidth}
		\includegraphics[width=\textwidth]{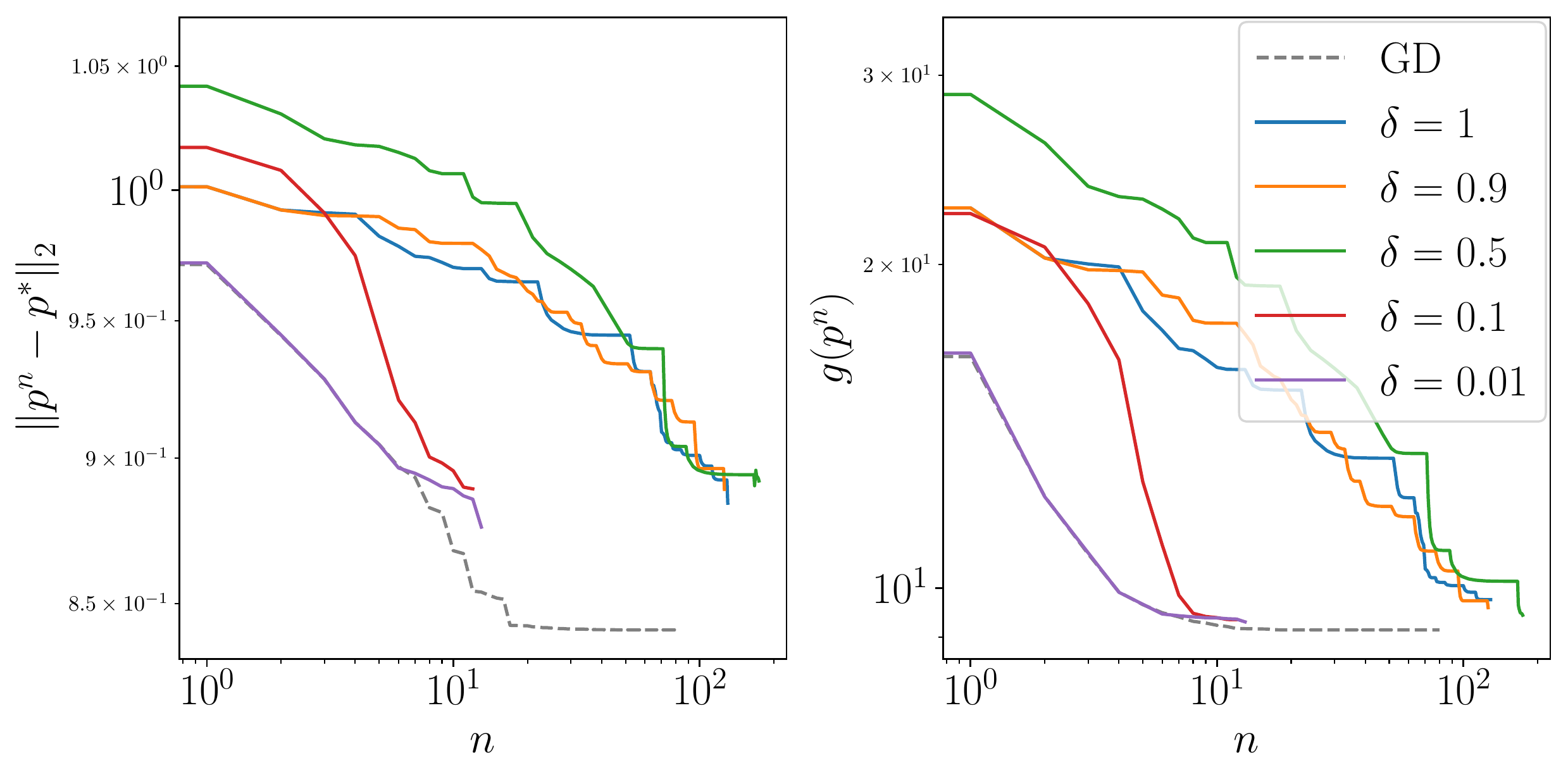}
		\caption{$N=4$; $\textup{dim}(P)=16$}
	\end{subfigure}

	\begin{subfigure}{0.7\textwidth}
		\includegraphics[width=\textwidth]{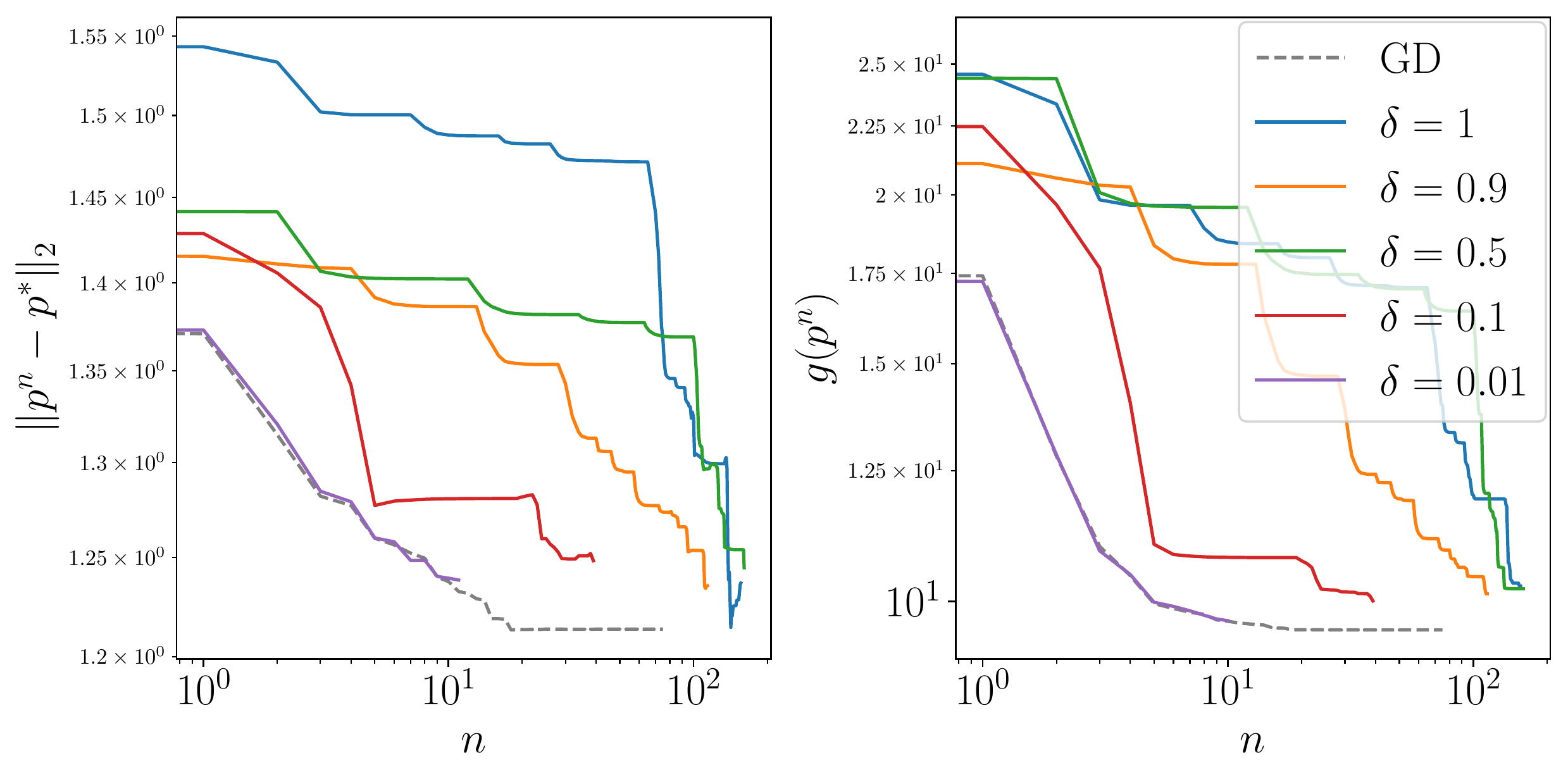}
		\caption{$N=8$; $\textup{dim}(P)=64$} 
	\end{subfigure}

	\caption{Results for the groundwater flow example from \cref{sec:gwf}, for a variety of parameter dimensions.
	In each row, the left-hand plot shows the distance between the parameter found at iteration $n$ and the true value $p^*$ of the parameter.
	The right-hand plot shows the value of the objective function, the negative log-likelihood in the Bayesian inference problem.} \label{fig:gwf_results}
\end{figure}

\begin{figure}
    \centering
    \includegraphics[width=\textwidth]{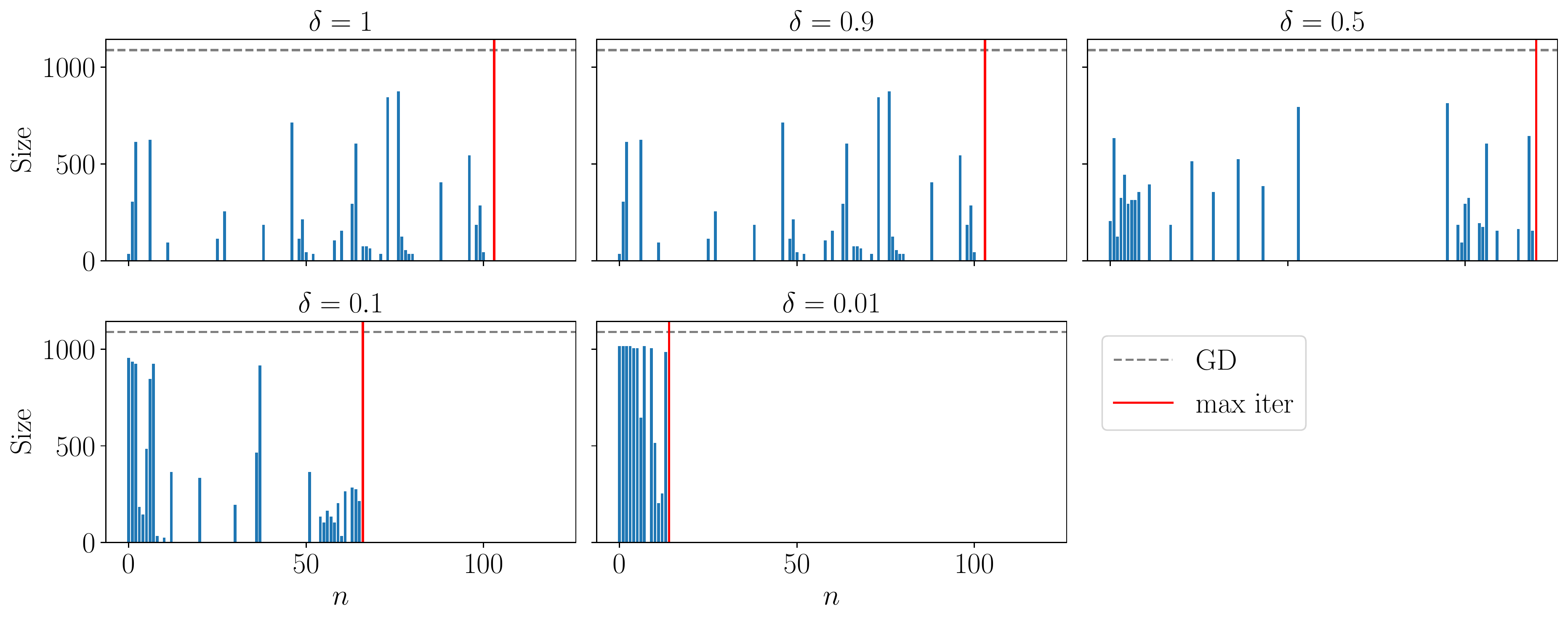}
    \caption{
    Size of the matrix whose Cholesky factorisation must be computed for each iteration of gradient descent in the groundwater flow example from \cref{sec:gwf}, as the parameter $\delta$ is varied.
    The figure is for $N=2$, $\textup{dim}(P)=4$, but results for other parameter dimensions are similar.
    The dashed gray line shows the size of the (sparse) problem that must be solved for the classical approach, while the red line indicates the iteration number at which convergence was achieved.
    }
    \label{fig:gwf_inversion_size}
\end{figure}

The results of the optimisation are displayed in \cref{fig:gwf_results}.
As in \cref{sec:fitzhugh_nagumo} one can clearly see the behaviour of the method reverting to that of gradient descent as the size of $\delta$ is decreased.
Further, performance appears to be broadly similar as the parameter dimension increases, reflecting that only a single function $\beta(x)$ must be learned, rather than $\frac{\wrt U}{\wrt p_i}$ for $i=1,\dots,\textup{dim}(P)$ as would be required in the forward approach.
Thus, the output dimension of the inferred function is independent of the parameter dimension.
While the \emph{input} dimension does grow with $\textup{dim}(P)$, for the purposes of the gradient descent algorithm, at iteration $n$ only the quality of inferences at and in the region of $p^n$ is relevant.
Since these points concentrate near $p^*$, performance does not appear to decay as the input dimension grows.

\cref{fig:gwf_inversion_size} compares the cost of the probabilistic approach with that of the classical approach, for $N=2$, $\textup{dim}(P) = 4$, by plotting the size of the matrix whose Cholesky factorisation that must be computed at each iteration in order to update the Cholesky factorisation of the Gram matrix with novel information, as discussed in \cref{sec:algorithm_discussion}.
We note that in general more information seems to be required than for the Fitzhugh-Nagumo example, so that the $10,000 \times 10,000$ limit on the size of the Gram matrix discussed in \cref{sec:algorithm_discussion} is generally what causes the algorithm to terminate, though from \cref{fig:gwf_results} it is clear that nevertheless \textsc{probjac} is \emph{close} to convergence when this occurs.
The higher cost is perhaps due to the fact that the right-hand side, $\frac{\partial g}{\partial u}$, is highly localised in this example.
It is nevertheless the case that throughout the gradient descent procedure, the size of the inversion problem that must be computed with the probabilistic approach is significantly smaller than that which must be computed with the classical approach, though since the matrix inverted in the classical approach is sparse the costs are not directly comparable.
Furthermore, as in \cref{sec:fitzhugh_nagumo}, for larger values of $\delta$ the approach shows the behaviour of being able to conduct a large number of iterations without needing to collect \emph{any} evaluations of the right-hand-side, due to the fact that the model is global over parameter space.

\section{Conclusion} \label{sec:conclusion}

In this paper we have presented a probabilistic approach to computing local sensitivities of differential equation models in both the forward and adjoint modes.
We presented an approach for incorporating these probabilistic gradients into a gradient descent algorithm, and examined the properties of this algorithm on two challenging applied problems with favourable results compared to classical approaches.
The chief advantages of the approach are that (i) gradients can be calculated at a lower cost than in classical approaches, (ii) that a global model for the gradient across parameter space is constructed, allowing for re-use of computational effort from previous iterations of gradient descent and (iii) that a full probability model is output, providing an error indicator that we used both to determine when to refine the approximation and to perform line searches.

Several possible avenues for future work present themselves.
The first would be continuing to develop applications of this algorithm within optimisation, either by developing versions of more sophisticated gradient-based optimisation algorithms which exploit probabilistic gradients, or by extending the framework to obtain higher order information to accelerate the optimisation.
Another would be to explore the use of probabilistic gradients in other applications.
In particular, we note that while computing the MAP point is an important problem in Bayesian inference, sophisticated Markov-chain Monte-Carlo algorithms for sampling the posterior also make use of this information, and the posterior distribution over the gradient presented herein could straightforwardly be incorporated into such algorithms.

\bibliographystyle{\thebibstyle}
\bibliography{refs}

\begin{thebibliography}{55}
\providecommand{\natexlab}[1]{#1}
\providecommand{\url}[1]{\texttt{#1}}
\expandafter\ifx\csname urlstyle\endcsname\relax
  \providecommand{\doi}[1]{doi: #1}\else
  \providecommand{\doi}{doi: \begingroup \urlstyle{rm}\Url}\fi

\bibitem[{\'{A}}lvarez et~al.(2012){\'{A}}lvarez, Rosasco, and
  Lawrence]{Alvarez2012}
M.~A. {\'{A}}lvarez, L.~Rosasco, and N.~D. Lawrence.
\newblock Kernels for vector-valued functions: a review.
\newblock \emph{Foundations and Trends in Machine Learning}, 4\penalty0
  (3):\penalty0 195--266, 2012.
\newblock \doi{10.1561/2200000036}.

\bibitem[Arriola and Hyman(2009)]{arriola2009sensitivity}
L.~Arriola and J.~M. Hyman.
\newblock Sensitivity analysis for uncertainty quantification in mathematical
  models.
\newblock In \emph{Mathematical and statistical estimation approaches in
  epidemiology}, pages 195--247. Springer, 2009.

\bibitem[Beddows et~al.(2017)Beddows, Kitwiroon, Williams, and
  Beevers]{Beddows2017}
A.~V. Beddows, N.~Kitwiroon, M.~L. Williams, and S.~D. Beevers.
\newblock Emulation and sensitivity analysis of the community multiscale air
  quality model for a {UK} ozone pollution episode.
\newblock \emph{Environmental Science {\&} Technology}, 51\penalty0
  (11):\penalty0 6229--6236, May 2017.
\newblock \doi{10.1021/acs.est.6b05873}.

\bibitem[Benner et~al.(2014)Benner, Sachs, and Volkwein]{benner2014model}
P.~Benner, E.~Sachs, and S.~Volkwein.
\newblock Model order reduction for {PDE} constrained optimization.
\newblock In \emph{Trends in PDE constrained optimization}, pages 303--326.
  Springer, 2014.

\bibitem[Benner et~al.(2015)Benner, Gugercin, and Willcox]{benner2015survey}
P.~Benner, S.~Gugercin, and K.~Willcox.
\newblock A survey of projection-based model reduction methods for parametric
  dynamical systems.
\newblock \emph{SIAM review}, 57\penalty0 (4):\penalty0 483--531, 2015.

\bibitem[Berlinet and Thomas-Agnan(2004)]{Berlinet2004}
A.~Berlinet and C.~Thomas-Agnan.
\newblock \emph{Reproducing Kernel Hilbert Spaces in Probability and
  Statistics}.
\newblock Springer {US}, 2004.
\newblock \doi{10.1007/978-1-4419-9096-9}.

\bibitem[Biegler et~al.(2007)Biegler, Ghattas, Heinkenschloss, Keyes, and van
  Bloemen~Waanders]{biegler2007real}
L.~T. Biegler, O.~Ghattas, M.~Heinkenschloss, D.~Keyes, and B.~van
  Bloemen~Waanders.
\newblock \emph{Real-time PDE-constrained Optimization}.
\newblock SIAM, 2007.

\bibitem[Bogachev(1998)]{Bogachev:1998}
V.~I. Bogachev.
\newblock \emph{Gaussian Measures}, volume~62.
\newblock American Mathematical Society Providence, 1998.

\bibitem[Bonnans and Shapiro(2013)]{bonnans2013perturbation}
J.~F. Bonnans and A.~Shapiro.
\newblock \emph{Perturbation analysis of optimization problems}.
\newblock Springer Science \& Business Media, 2013.

\bibitem[Cheng et~al.(2020)Cheng, Lu, Ling, and Zhou]{Cheng2020}
K.~Cheng, Z.~Lu, C.~Ling, and S.~Zhou.
\newblock Surrogate-assisted global sensitivity analysis: an overview.
\newblock \emph{Structural and Multidisciplinary Optimization}, 61\penalty0
  (3):\penalty0 1187--1213, Jan. 2020.
\newblock \doi{10.1007/s00158-019-02413-5}.

\bibitem[Cialenco et~al.(2012)Cialenco, Fasshauer, and Ye]{Cialenco2012}
I.~Cialenco, G.~E. Fasshauer, and Q.~Ye.
\newblock Approximation of stochastic partial differential equations by a
  kernel-based collocation method.
\newblock \emph{Int. J. Comput. Math.}, 89\penalty0 (18):\penalty0 2543--2561,
  2012.
\newblock \doi{10.1080/00207160.2012.688111}.

\bibitem[Cleary et~al.(2020)Cleary, Garbuno-Inigo, Lan, Schneider, and
  Stuart]{CLEARY2020109716}
E.~Cleary, A.~Garbuno-Inigo, S.~Lan, T.~Schneider, and A.~M. Stuart.
\newblock Calibrate, {Emulate}, {Sample}.
\newblock \emph{Journal of Computational Physics}, page 109716, 2020.
\newblock ISSN 0021-9991.
\newblock \doi{https://doi.org/10.1016/j.jcp.2020.109716}.

\bibitem[Cockayne(2019)]{CockayneThesis}
J.~Cockayne.
\newblock \emph{Bayesian Probabilistic Numerical Methods}.
\newblock PhD thesis, University of Warwick, 2019.

\bibitem[Cockayne et~al.(2019)Cockayne, Oates, Sullivan, and
  Girolami]{Cockayne2019}
J.~Cockayne, C.~J. Oates, T.~J. Sullivan, and M.~Girolami.
\newblock Bayesian probabilistic numerical methods.
\newblock \emph{{SIAM} Review}, 61\penalty0 (3):\penalty0 756--789, Jan. 2019.
\newblock \doi{10.1137/17m1139357}.

\bibitem[Curry(1944)]{Curry1944}
H.~B. Curry.
\newblock The method of steepest descent for non-linear minimization problems.
\newblock \emph{Q APPL MATH}, 2\penalty0 (3):\penalty0 258--261, Oct. 1944.
\newblock \doi{10.1090/qam/10667}.

\bibitem[Drohmann and Carlberg(2015)]{drohmann2015romes}
M.~Drohmann and K.~Carlberg.
\newblock The {ROMES} method for statistical modeling of reduced-order-model
  error.
\newblock \emph{SIAM/ASA Journal on Uncertainty Quantification}, 3\penalty0
  (1):\penalty0 116--145, 2015.

\bibitem[Evans(2010)]{Evans2010}
L.~Evans.
\newblock \emph{Partial Differential Equations}.
\newblock American Mathematical Society, Mar. 2010.
\newblock \doi{10.1090/gsm/019}.

\bibitem[Fisher et~al.(2009)Fisher, Nocedal, Tr{\'e}molet, and
  Wright]{fisher2009data}
M.~Fisher, J.~Nocedal, Y.~Tr{\'e}molet, and S.~J. Wright.
\newblock Data assimilation in weather forecasting: a case study in
  {PDE}-constrained optimization.
\newblock \emph{Optimization and Engineering}, 10\penalty0 (3):\penalty0
  409--426, 2009.

\bibitem[FitzHugh(1961)]{FitzHugh1961}
R.~FitzHugh.
\newblock Impulses and physiological states in theoretical models of nerve
  membrane.
\newblock \emph{Biophysical Journal}, 1\penalty0 (6):\penalty0 445--466, July
  1961.
\newblock \doi{10.1016/s0006-3495(61)86902-6}.

\bibitem[Geyer(2011)]{geyer2011introduction}
C.~Geyer.
\newblock Introduction to markov chain monte carlo.
\newblock \emph{Handbook of markov chain monte carlo}, 20116022:\penalty0 45,
  2011.

\bibitem[Girard et~al.(2016)Girard, Mallet, Korsakissok, and
  Mathieu]{Girard2016}
S.~Girard, V.~Mallet, I.~Korsakissok, and A.~Mathieu.
\newblock Emulation and {Sobol$'$} sensitivity analysis of an atmospheric
  dispersion model applied to the {Fukushima} nuclear accident.
\newblock \emph{Journal of Geophysical Research: Atmospheres}, 121\penalty0
  (7):\penalty0 3484--3496, Apr. 2016.
\newblock \doi{10.1002/2015jd023993}.

\bibitem[Gunzburger(2002)]{gunzburgerperspective}
M.~Gunzburger.
\newblock \emph{Perspective in flow control and optimization (2003)}.
\newblock SIAM, Philadelphia, 2002.

\bibitem[Hairer et~al.(1993)Hairer, Nørsett, and Wanner]{Hairer1993}
E.~Hairer, S.~Nørsett, and G.~Wanner.
\newblock \emph{Solving Ordinary Differential Equations {I}: Nonstiff
  Problems}.
\newblock Springer, 1993.

\bibitem[Hartman and Mestha(2017)]{hartman2017deep}
D.~Hartman and L.~K. Mestha.
\newblock A deep learning framework for model reduction of dynamical systems.
\newblock In \emph{2017 IEEE Conference on Control Technology and Applications
  (CCTA)}, pages 1917--1922. IEEE, 2017.

\bibitem[Hennig et~al.(2015)Hennig, Osborne, and Girolami]{Hennig2015a}
P.~Hennig, M.~A. Osborne, and M.~Girolami.
\newblock Probabilistic numerics and uncertainty in computations.
\newblock \emph{J. R. Stat. Soc. A Stat.}, 471\penalty0 (2179):\penalty0
  20150142, 17, 2015.
\newblock \doi{10.1098/rspa.2015.0142}.

\bibitem[Herzog and Kunisch(2010)]{Herzog2010}
R.~Herzog and K.~Kunisch.
\newblock Algorithms for {PDE}-constrained optimization.
\newblock \emph{{GAMM}-Mitteilungen}, 33\penalty0 (2):\penalty0 163--176, Oct.
  2010.
\newblock \doi{10.1002/gamm.201010013}.

\bibitem[Higdon et~al.(2004)Higdon, Kennedy, Cavendish, Cafeo, and
  Ryne]{higdon2004combining}
D.~Higdon, M.~Kennedy, J.~C. Cavendish, J.~A. Cafeo, and R.~D. Ryne.
\newblock Combining field data and computer simulations for calibration and
  prediction.
\newblock \emph{SIAM Journal on Scientific Computing}, 26\penalty0
  (2):\penalty0 448--466, 2004.

\bibitem[Ito and Kunisch(2008)]{ito2008lagrange}
K.~Ito and K.~Kunisch.
\newblock \emph{Lagrange multiplier approach to variational problems and
  applications}.
\newblock SIAM, 2008.

\bibitem[Jin et~al.(2004)Jin, Chen, and Sudjianto]{Jin2004}
R.~Jin, W.~Chen, and A.~Sudjianto.
\newblock Analytical metamodel-based global sensitivity analysis and
  uncertainty propagation for robust design.
\newblock In \emph{{SAE} Technical Paper Series}. {SAE} International, Mar.
  2004.
\newblock \doi{10.4271/2004-01-0429}.

\bibitem[Kennedy and O'Hagan(2001)]{Kennedy2001}
M.~C. Kennedy and A.~O'Hagan.
\newblock Bayesian calibration of computer models.
\newblock \emph{Journal of the Royal Statistical Society: Series B (Statistical
  Methodology)}, 63\penalty0 (3):\penalty0 425--464, Aug. 2001.
\newblock \doi{10.1111/1467-9868.00294}.

\bibitem[Lan et~al.(2016)Lan, Bui-Thanh, Christie, and
  Girolami]{lan2016emulation}
S.~Lan, T.~Bui-Thanh, M.~Christie, and M.~Girolami.
\newblock Emulation of higher-order tensors in manifold {Monte Carlo} methods
  for {Bayesian} inverse problems.
\newblock \emph{Journal of Computational Physics}, 308:\penalty0 81--101, 2016.

\bibitem[Mahsereci and Hennig(2015)]{Mahsereci2015}
M.~Mahsereci and P.~Hennig.
\newblock Probabilistic line searches for stochastic optimization.
\newblock In C.~Cortes, N.~D. Lawrence, D.~D. Lee, M.~Sugiyama, and R.~Garnett,
  editors, \emph{Advances in Neural Information Processing Systems 28}, pages
  181--189. Curran Associates, Inc., 2015.
\newblock URL
  \url{http://papers.nips.cc/paper/5753-probabilistic-line-searches-for-stochastic-optimization.pdf}.

\bibitem[Newman~III et~al.(1999)Newman~III, Taylor~III, Barnwell, Newman, and
  Hou]{newman1999overview}
J.~C. Newman~III, A.~C. Taylor~III, R.~W. Barnwell, P.~A. Newman, and G.~J.-W.
  Hou.
\newblock Overview of sensitivity analysis and shape optimization for complex
  aerodynamic configurations.
\newblock \emph{Journal of Aircraft}, 36\penalty0 (1):\penalty0 87--96, 1999.

\bibitem[Nocedal and Wright(2006)]{Nocedal2006}
J.~Nocedal and S.~J. Wright.
\newblock \emph{Numerical Optimization}.
\newblock Springer New York, 2006.
\newblock \doi{10.1007/978-0-387-40065-5}.
\newblock URL \url{https://doi.org/10.1007/978-0-387-40065-5}.

\bibitem[Oakley and O'Hagan(2002)]{oakley2002bayesian}
J.~Oakley and A.~O'Hagan.
\newblock Bayesian inference for the uncertainty distribution of computer model
  outputs.
\newblock \emph{Biometrika}, 89\penalty0 (4):\penalty0 769--784, 2002.

\bibitem[Oakley and O'Hagan(2004)]{Oakley2004}
J.~E. Oakley and A.~O'Hagan.
\newblock Probabilistic sensitivity analysis of complex models: a {Bayesian}
  approach.
\newblock \emph{Journal of the Royal Statistical Society: Series B (Statistical
  Methodology)}, 66\penalty0 (3):\penalty0 751--769, Aug. 2004.
\newblock \doi{10.1111/j.1467-9868.2004.05304.x}.

\bibitem[Oates and Sullivan(2019)]{Oates2019}
C.~J. Oates and T.~J. Sullivan.
\newblock A modern retrospective on probabilistic numerics.
\newblock \emph{Statistics and Computing}, 29\penalty0 (6):\penalty0
  1335--1351, Oct. 2019.
\newblock \doi{10.1007/s11222-019-09902-z}.

\bibitem[Osborne(2010)]{Osborne:2010}
M.~Osborne.
\newblock \emph{Bayesian {Gaussian} {Processes} for {Sequential} {Prediction},
  {Optimisation} and {Quadrature}}.
\newblock PhD thesis, PhD thesis, University of Oxford, 2010.

\bibitem[Owhadi(2015)]{Owhadi2015}
H.~Owhadi.
\newblock Bayesian numerical homogenization.
\newblock \emph{Multiscale Modeling {\&} Simulation}, 13\penalty0 (3):\penalty0
  812--828, Jan. 2015.
\newblock \doi{10.1137/140974596}.

\bibitem[Owhadi(2017)]{Owhadi2017a}
H.~Owhadi.
\newblock Multigrid with rough coefficients and multiresolution operator
  decomposition from hierarchical information games.
\newblock \emph{{SIAM} Review}, 59\penalty0 (1):\penalty0 99--149, Jan. 2017.
\newblock \doi{10.1137/15m1013894}.

\bibitem[Owhadi and Zhang(2017)]{Owhadi2017b}
H.~Owhadi and L.~Zhang.
\newblock Gamblets for opening the complexity-bottleneck of implicit schemes
  for hyperbolic and parabolic {ODEs}/{PDEs} with rough coefficients.
\newblock \emph{Journal of Computational Physics}, 347:\penalty0 99--128, Oct.
  2017.
\newblock \doi{10.1016/j.jcp.2017.06.037}.

\bibitem[Pulch et~al.(2015)Pulch, ter Maten, and
  Augustin]{pulch2015sensitivity}
R.~Pulch, E.~J.~W. ter Maten, and F.~Augustin.
\newblock Sensitivity analysis and model order reduction for random linear
  dynamical systems.
\newblock \emph{Mathematics and Computers in Simulation}, 111:\penalty0 80--95,
  2015.

\bibitem[Rasmussen and Williams(2005)]{Rasmussen2005}
C.~E. Rasmussen and C.~K.~I. Williams.
\newblock \emph{Gaussian Processes for Machine Learning}.
\newblock The {MIT} Press, 2005.
\newblock \doi{10.7551/mitpress/3206.001.0001}.

\bibitem[Renardy et~al.(2018)Renardy, Yi, Xiu, and Chou]{Renardy2018}
M.~Renardy, T.-M. Yi, D.~Xiu, and C.-S. Chou.
\newblock Parameter uncertainty quantification using surrogate models applied
  to a spatial model of yeast mating polarization.
\newblock \emph{{PLOS} Computational Biology}, 14\penalty0 (5):\penalty0
  e1006181, May 2018.
\newblock \doi{10.1371/journal.pcbi.1006181}.

\bibitem[Sacks et~al.(1989)Sacks, Welch, Mitchell, and Wynn]{sacks1989design}
J.~Sacks, W.~J. Welch, T.~J. Mitchell, and H.~P. Wynn.
\newblock Design and analysis of computer experiments.
\newblock \emph{Statistical science}, pages 409--423, 1989.

\bibitem[San and Maulik(2018{\natexlab{a}})]{san2018extreme}
O.~San and R.~Maulik.
\newblock Extreme learning machine for reduced order modeling of turbulent
  geophysical flows.
\newblock \emph{Physical Review E}, 97\penalty0 (4):\penalty0 042322,
  2018{\natexlab{a}}.

\bibitem[San and Maulik(2018{\natexlab{b}})]{san2018neural}
O.~San and R.~Maulik.
\newblock Neural network closures for nonlinear model order reduction.
\newblock \emph{Advances in Computational Mathematics}, 44\penalty0
  (6):\penalty0 1717--1750, 2018{\natexlab{b}}.

\bibitem[Sengupta et~al.(2014)Sengupta, Friston, and
  Penny]{sengupta2014efficient}
B.~Sengupta, K.~J. Friston, and W.~D. Penny.
\newblock Efficient gradient computation for dynamical models.
\newblock \emph{NeuroImage}, 98:\penalty0 521--527, 2014.

\bibitem[Shafto et~al.(2012)Shafto, Conroy, Doyle, Glaessgen, Kemp, LeMoigne,
  and Wang]{shafto2012modeling}
M.~Shafto, M.~Conroy, R.~Doyle, E.~Glaessgen, C.~Kemp, J.~LeMoigne, and
  L.~Wang.
\newblock Modeling, simulation, information technology \& processing roadmap.
\newblock \emph{National Aeronautics and Space Administration}, 2012.

\bibitem[Sheriffdeen et~al.(2019)Sheriffdeen, Ragusa, Morel, Adams, and
  Bui-Thanh]{sheriffdeen2019accelerating}
S.~Sheriffdeen, J.~C. Ragusa, J.~E. Morel, M.~L. Adams, and T.~Bui-Thanh.
\newblock Accelerating {PDE}-constrained inverse solutions with {Deep Learning}
  and {Reduced Order Models}.
\newblock \emph{arXiv preprint arXiv:1912.08864}, 2019.

\bibitem[Sobol$'$(2001)]{Sobol2001}
I.~Sobol$'$.
\newblock Global sensitivity indices for nonlinear mathematical models and
  their {Monte Carlo} estimates.
\newblock \emph{Mathematics and Computers in Simulation}, 55\penalty0
  (1-3):\penalty0 271--280, Feb. 2001.
\newblock \doi{10.1016/s0378-4754(00)00270-6}.

\bibitem[Stuart and Teckentrup(2018)]{stuart2018posterior}
A.~Stuart and A.~Teckentrup.
\newblock Posterior consistency for {Gaussian} process approximations of
  {Bayesian} posterior distributions.
\newblock \emph{Mathematics of Computation}, 87\penalty0 (310):\penalty0
  721--753, 2018.

\bibitem[Stuart(2010)]{Stuart:2010}
A.~M. Stuart.
\newblock Inverse problems: A {Bayesian} perspective.
\newblock \emph{Acta Numerica}, 19:\penalty0 451--559, may 2010.
\newblock \doi{10.1017/s0962492910000061}.

\bibitem[Villani(2009)]{Villani2009}
C.~Villani.
\newblock \emph{Optimal Transport}.
\newblock Springer Berlin Heidelberg, 2009.
\newblock \doi{10.1007/978-3-540-71050-9}.

\bibitem[Wendland(2004)]{Wendland2004}
H.~Wendland.
\newblock \emph{Scattered Data Approximation}.
\newblock Cambridge University Press, Dec. 2004.
\newblock \doi{10.1017/cbo9780511617539}.

\end{thebibliography}


\begin{thebibliography}{6}
\providecommand{\natexlab}[1]{#1}
\providecommand{\url}[1]{\texttt{#1}}
\expandafter\ifx\csname urlstyle\endcsname\relax
  \providecommand{\doi}[1]{doi: #1}\else
  \providecommand{\doi}{doi: \begingroup \urlstyle{rm}\Url}\fi

\bibitem[Gelbrich(1990)]{Gelbrich1990}
M.~Gelbrich.
\newblock On a formula for the {L2} {Wasserstein} metric between measures on
  {Euclidean} and {Hilbert} spaces.
\newblock \emph{Mathematische Nachrichten}, 147\penalty0 (1):\penalty0
  185--203, 1990.
\newblock \doi{10.1002/mana.19901470121}.
\newblock URL \url{https://doi.org/10.1002/mana.19901470121}.

\bibitem[Gohberg et~al.(2003)Gohberg, Goldberg, and Kaashoek]{Gohberg2003}
I.~Gohberg, S.~Goldberg, and M.~A. Kaashoek.
\newblock \emph{Basic Classes of Linear Operators}.
\newblock Birkh\"{a}user Basel, 2003.
\newblock \doi{10.1007/978-3-0348-7980-4}.

\bibitem[Masarotto et~al.(2018)Masarotto, Panaretos, and Zemel]{Masarotto2018}
V.~Masarotto, V.~M. Panaretos, and Y.~Zemel.
\newblock Procrustes metrics on covariance operators and optimal transportation
  of {Gaussian} processes.
\newblock \emph{Sankhya A}, 81\penalty0 (1):\penalty0 172--213, May 2018.
\newblock \doi{10.1007/s13171-018-0130-1}.

\bibitem[Nocedal and Wright(2006)]{Nocedal2006}
J.~Nocedal and S.~J. Wright.
\newblock \emph{Numerical Optimization}.
\newblock Springer New York, 2006.
\newblock \doi{10.1007/978-0-387-40065-5}.
\newblock URL \url{https://doi.org/10.1007/978-0-387-40065-5}.

\bibitem[Osborne(2010)]{Osborne:2010}
M.~Osborne.
\newblock \emph{Bayesian {Gaussian} {Processes} for {Sequential} {Prediction},
  {Optimisation} and {Quadrature}}.
\newblock PhD thesis, PhD thesis, University of Oxford, 2010.

\bibitem[Pigoli et~al.(2014)Pigoli, Aston, Dryden, and Secchi]{Pigoli2014}
D.~Pigoli, J.~A.~D. Aston, I.~L. Dryden, and P.~Secchi.
\newblock Distances and inference for covariance operators.
\newblock \emph{Biometrika}, 101\penalty0 (2):\penalty0 409--422, Apr. 2014.
\newblock \doi{10.1093/biomet/asu008}.

\end{thebibliography}

\ifarxiv


\section*{Supplementary material}

# Supplementary Material

For: Probabilistic Gradients for Fast Calibration of Differential Equation Models.

Jon Cockayne*    Andrew B. Duncan†

February 22, 2021

## S1 Proofs

In this section several important lemmas are presented in Section S1.1, then the proofs of the main theoretical results are presented in Section S1.2.

### S1.1 Lemmas

First we present several Lemmas that will be of use in the proofs of the main theoretical results. The first Lemma provides a representation of the conditional mean and covariance of Gaussian measures that will be utilised in the main proofs.

**Lemma S1.** *Let $\bar\mu \sim \mathcal{N}(\bar a, \bar C)$ be a conditional Gaussian measure on $\mathcal{H}$, with mean and covariance given by*

$$\bar a = a + C\mathcal{I}^\dagger [\mathcal{I} C \mathcal{I}^\dagger]^{-1}(f - \mathcal{I}a)$$
$$\bar C = C - C\mathcal{I}^\dagger [\mathcal{I} C \mathcal{I}^\dagger]^{-1} \mathcal{I} C.$$

*Here $a \in \mathcal{H}$, $C$ is a covariance operator on $\mathcal{H}$ with associated RKHS $\mathcal{H}_C$, $\mathcal{I} : \mathcal{H}_C \to \mathbb{R}^d$ is a bounded linear operator and $f = \mathcal{I}u^\dagger$ for some $u^\dagger \in \mathcal{H}_C$. Let*

$$P = C\mathcal{I}^\dagger [\mathcal{I} C \mathcal{I}^\dagger]^{-1} \mathcal{I}$$

*Then we have the following*

1. *$P$ is an orthogonal projection in the inner product induced by $C^{-1}$.*

2. *$\bar a = Pu^\dagger + (\text{Id} - P)a$*

3. *$\bar C = (\text{Id} - P)C$*

---
*The Alan Turing Institute (jcockayne@turing.ac.uk)
†Imperial College London and The Alan Turing Institute (a.duncan@imperial.ac.uk)

4. $\bar{C} = (\text{Id} - P)C(\text{Id} - P)^\dagger$

*Proof.* For (1.), it is straightforward to verify that $P^2 = P$, so that $P$ is a projection. To verify that $P$ is an *orthogonal* projection recall that a projection is orthogonal if and only if it is self-adjoint. To verify this we have

$$\langle Pu, v \rangle_{C^{-1}} = \left\langle C^{-1} C \mathcal{I}^\dagger [\mathcal{I} C \mathcal{I}^\dagger]^{-1} \mathcal{I} u, v \right\rangle$$
$$= \left\langle u, \mathcal{I}^\dagger [\mathcal{I} C \mathcal{I}^\dagger]^{-1} \mathcal{I} v \right\rangle$$
$$= \langle u, Pv \rangle_{C^{-1}}.$$

(2.) and (3.) are clear by inspection. For (4.) note that

$$(\text{Id} - P)C(\text{Id} - P)^\dagger = C - PC - CP^\dagger + PCP^\dagger$$
$$PC = C\mathcal{I}^\dagger [\mathcal{I} C \mathcal{I}^\dagger]^{-1} \mathcal{I} C$$
$$CP^\dagger = C\mathcal{I}^\dagger [\mathcal{I} C \mathcal{I}^\dagger]^{-1} \mathcal{I} = PC$$
$$\implies PCP^\dagger = P^2 C = PC$$

So that $C - PC - CP^\dagger + PCP^\dagger = C - PC = \bar{C}$, from (3.). □

The next lemma provides a bound for the Wasserstein distance that will be used in Theorem 3.4.

**Lemma S2.** *Let $\mu_1 = \mathcal{N}(a_1, C_1)$ and $\mu_2 = \mathcal{N}(a_2, C_2)$ each be Gaussian measures on $\mathcal{U}$. Let $S_1, S_2$ each be square roots of $C_1, C_2$ such that $C_1 = S_1^\dagger S_1$ and $C_2 = S_2^\dagger S_2$ Then we have that*

$$W_2^2(\mu_1, \mu_2) \leq \|a_1 - a_2\|_{\mathcal{H}}^2 + \|S_1 - S_2\|_{\mathsf{HS}}^2.$$

*Proof.* Let $\mu_1' = \mathcal{N}(0, C_1)$ and $\mu_2' = \mathcal{N}(0, C_2)$. From Gelbrich [1990] we have that

$$W_2(\mu_1, \mu_2)^2 = \|a_1 - a_2\|^2 + W_2(\mu_1', \mu_2')^2.$$

Now we have that for *centered* Gaussian measures $\mu_1 = \mathcal{N}(0, C_1)$, $\mu_2 = \mathcal{N}(0, C_2)$

$$W_2(\mu_1, \mu_2) = \inf_{R: R^\dagger R = \text{Id}} \|S_1 - RS_2\|_{\mathsf{HS}}. \tag{S1}$$

where the fact that this holds for square roots other than the *symmetric* square roots is clear from the fact that the Wasserstein distance above is equal to the Procrustes distance between their respective covariance operators [Masarotto et al., 2018, Proposition 3], and the Procrustes distance need not be expressed in terms of symmetric square roots [Pigoli et al., 2014]. By taking $R = \text{Id}$ we obtain the following bound:

$$W_2(\mu_1', \mu_2')^2 \leq \|S_1 - S_2\|_{\mathsf{HS}}^2.$$

□

## S1.2 Proofs of Main Results

*Proof of Theorem 3.3.* From Theorem S1 we have

$$\begin{aligned}
|\mathcal{L}\bar{a} - \mathcal{L}u^*| &= |\mathcal{L}([\text{Id} - P]a + Pu^* - u^*)| \\
&= |\mathcal{L}(\text{Id} - P)(a - u^*)| \\
&= \langle (\text{Id} - P)(a - u^*), \ell \rangle \\
&= \langle (\text{Id} - P)(a - u^*), C\ell \rangle_{C^{-1}} \\
&= \langle a - u^*, \bar{C}\ell \rangle_{C^{-1}} \\
&\leq \|\bar{C}\ell\|_{C^{-1}} \|a - u^*\|_{C^{-1}}.
\end{aligned}$$

where $\ell$ is the Riesz representer of $\mathcal{L}$. On the fourth line the fact that the projection $P$ is self-adjoint in the $C^{-1}$-inner product was used, and the final line is from the Cauchy–Schwarz inequality. Now, note that by substituting back in the expression for $\bar{C}$ we have

$$\begin{aligned}
\|\bar{C}\ell\|^2_{C^{-1}} &= \langle \bar{C}\ell, \bar{C}\ell \rangle_{C^{-1}} \\
&= \langle C\ell, C\ell \rangle_{C^{-1}} - 2\langle PC\ell, C\ell \rangle_{C^{-1}} + \langle PC\ell, PC\ell \rangle_{C^{-1}} \\
&= \langle C\ell, \ell \rangle - \langle PC\ell, \ell \rangle \\
&= \langle (\text{Id} - P)C\ell, \ell \rangle \\
&= \langle \bar{C}\ell, \ell \rangle = \mathcal{L}\bar{C}\ell = \mathcal{L}\bar{C}\mathcal{L}^\dagger
\end{aligned}$$

where on the third line we have used the fact that $P$ is an orthogonal projection, and the final line comes from interpreting $\mathcal{L}\bar{C}\ell$ as an operator from $\mathbb{R} \to \mathbb{R}$. This completes the proof. $\square$

*Proof of Theorem 3.4.* Let $G = \mathcal{I}C\mathcal{I}^\dagger$ and $\hat{G} = \hat{\mathcal{I}}C\hat{\mathcal{I}}^\dagger$. From Theorem S1 there exists an operator $P = C\mathcal{I}^\dagger G^{-1}\mathcal{I}$ such that

$$\begin{aligned}
\bar{a} &= Pu^\dagger + (\text{Id} - P)a \\
\bar{C} &= (\text{Id} - P)C(\text{Id} - P)^\dagger.
\end{aligned}$$

By the same logic, for $\hat{P} = C\hat{\mathcal{I}}^\dagger \hat{G}^{-1}\hat{\mathcal{I}}$ we have

$$\begin{aligned}
\hat{a} &= \hat{P}u^\dagger + (\text{Id} - \hat{P})a \\
\hat{C} &= (\text{Id} - \hat{P})C(\text{Id} - \hat{P})^\dagger.
\end{aligned}$$

Therefore, (assymmetric) square roots of both $\bar{C}$ and $\hat{C}$ are given by

$$\begin{aligned}
\bar{C}^{\frac{1}{2}} &= C^{\frac{1}{2}}(\text{Id} - \bar{P})^\dagger \\
\hat{C}^{\frac{1}{2}} &= C^{\frac{1}{2}}(\text{Id} - \hat{P})^\dagger
\end{aligned}$$

where $C^{\frac{1}{2}}$ is the symmetric[1] square root of $C$.

---

[1] Indeed, $C^{\frac{1}{2}}$ may also be an arbitrary square root of $C$, but the symmetric square root is used to simplify notation.

We begin by bounding the error in the posterior means. Recalling that $f = \mathcal{I} u^\dagger$ and $\hat{f} = \hat{\mathcal{I}} u^\dagger$, we have that

$$\bar{a} - \hat{a} = C\mathcal{I}^\dagger G^{-1} f - C\hat{\mathcal{I}}^\dagger \hat{G}^{-1} \hat{f} - (P - \hat{P})a$$
$$\implies \|\bar{a} - \hat{a}\|_{\mathcal{H}} \leq \|C(\mathcal{I}^\dagger G^{-1} f - \hat{\mathcal{I}}^\dagger \hat{G}^{-1} \hat{f})\|_{\mathbb{R}^d} + \|P - \hat{P}\|_{\mathcal{H}} \|a\|_{\mathcal{H}}.$$

Furthermore

$$\bar{C}^{\frac{1}{2}} - \hat{C}^{\frac{1}{2}} = C^{\frac{1}{2}}(P - \hat{P})$$
$$\implies \|\bar{C}^{\frac{1}{2}} - \hat{C}^{\frac{1}{2}}\|_{\mathsf{HS}} \leq \|P - \hat{P}\|_{\mathcal{H}} \|C^{\frac{1}{2}}\|_{\mathsf{HS}}$$

We thus have that

$$W_2(\bar{\mu}, \hat{\mu}) \leq \underbrace{\|C(\mathcal{I}^\dagger G^{-1} f - \hat{\mathcal{I}}^\dagger \hat{G}^{-1} \hat{f})\|_{\mathbb{R}^d}}_{(1)} + \underbrace{\|P - \hat{P}\|_{\mathcal{H}}}_{(2)}(\|a\|_{\mathcal{H}} + \|C^{\frac{1}{2}}\|_{\mathsf{HS}})$$

We first focus on (1). Clearly

$$\|C(\mathcal{I}^\dagger G^{-1} f - \hat{\mathcal{I}}^\dagger \hat{G}^{-1} \hat{f})\|_{\mathbb{R}^d} \leq \|C\|_{\mathcal{H}} \|\mathcal{I}^\dagger G^{-1} f - \hat{\mathcal{I}}^\dagger \hat{G}^{-1} \hat{f}\|_{\mathbb{R}^d}$$
$$\|\mathcal{I}^\dagger G^{-1} f - \hat{\mathcal{I}}^\dagger \hat{G}^{-1} \hat{f}\|_{\mathbb{R}^d} \leq \|\mathcal{I} - \hat{\mathcal{I}}\|_{\mathcal{H} \to \mathbb{R}^d} \|G^{-1} f\|_{\mathbb{R}^d}$$
$$+ \|\hat{\mathcal{I}}\|_{\mathcal{H} \to \mathbb{R}^d} \|G^{-1} - \hat{G}^{-1}\|_{\mathbb{R}^d} \|f\|_{\mathbb{R}^d}$$
$$+ \|\hat{\mathcal{I}}^\dagger \hat{G}^{-1}\|_{\mathbb{R}^d \to \mathcal{H}} \|f - \hat{f}\|_{\mathbb{R}^d}$$

We now focus on $\|P - \hat{P}\|_{\mathcal{H}}$. Clearly

$$\|P - \hat{P}\|_{\mathcal{H}} \leq \|\mathcal{I}^\dagger G^{-1} \mathcal{I} - \hat{\mathcal{I}}^\dagger \hat{G}^{-1} \hat{\mathcal{I}}\|_{\mathcal{H}} \|C\|_{\mathcal{H}}$$

and

$$\|\mathcal{I}^\dagger G^{-1} \mathcal{I} - \hat{\mathcal{I}}^\dagger \hat{G}^{-1} \hat{\mathcal{I}}\|_{\mathcal{H}} \leq \|\mathcal{I} - \hat{\mathcal{I}}\|_{\mathcal{H} \to \mathbb{R}^d} \|G^{-1} \mathcal{I}\|_{\mathcal{H}}$$
$$+ \|\mathcal{I}\|_{\mathcal{H} \to \mathbb{R}^d} \|\hat{\mathcal{I}}\|_{\mathcal{H} \to \mathbb{R}^d} \|G^{-1} - \hat{G}^{-1}\|_{\mathbb{R}^d}$$
$$+ \|\mathcal{I} - \hat{\mathcal{I}}\|_{\mathcal{H} \to \mathbb{R}^d} \|\hat{G}^{-1} \hat{\mathcal{I}}\|_{\mathcal{H} \to \mathbb{R}^d}.$$

It remains to bound $\|G^{-1} - \hat{G}^{-1}\|_{\mathbb{R}^d}$ in terms of $\|\mathcal{I} - \hat{\mathcal{I}}\|_{\mathcal{H} \to \mathbb{R}^d}$. This can be accomplished using Gohberg et al. [2003, Chapter II, Corollary 8.2], which states that whenever $\|G - \hat{G}\| \leq \frac{1}{\|G^{-1}\|}$, we may bound $\|G^{-1} - \hat{G}^{-1}\|$ as follows:

$$\|G^{-1} - \hat{G}^{-1}\|_{\mathbb{R}^d} \leq \frac{\|G^{-1}\|_{\mathbb{R}^d}^2 \|G - \hat{G}\|_{\mathbb{R}^d}}{1 - \|G^{-1}\|_{\mathbb{R}^d} \|G - \hat{G}\|_{\mathbb{R}^d}}$$

Recall that

$$\|G - \hat{G}\|_{\mathbb{R}^d} = \|\mathcal{I} C \mathcal{I}^\dagger - \hat{\mathcal{I}} C \hat{\mathcal{I}}^\dagger\|_{\mathbb{R}^d}$$
$$\leq \|(\mathcal{I} - \hat{\mathcal{I}}) C \mathcal{I}^\dagger\|_{\mathbb{R}^d} + \|\hat{\mathcal{I}} C (\mathcal{I} - \hat{\mathcal{I}})^\dagger\|_{\mathbb{R}^d}$$
$$\leq \|\mathcal{I} - \hat{\mathcal{I}}\|_{\mathcal{H} \to \mathbb{R}^d} (\|\mathcal{I} C\|_{\mathcal{H} \to \mathbb{R}^d} + \|\hat{\mathcal{I}} C\|_{\mathcal{H} \to \mathbb{R}^d})$$

so that, letting $\alpha = \|\mathcal{I}C\|_{\mathcal{H}\to\mathbb{R}^d} + \|\hat{\mathcal{I}}C\|_{\mathcal{H}\to\mathbb{R}^d}$,

$$\|G^{-1} - \hat{G}^{-1}\|_{\mathbb{R}^d} \leq \|G^{-1}\|_{\mathbb{R}^d} \frac{\alpha\|G^{-1}\|_{\mathbb{R}^d}\|\mathcal{I} - \hat{\mathcal{I}}\|_{\mathcal{H}\to\mathbb{R}^d}}{1 - \alpha\|G^{-1}\|_{\mathbb{R}^d}\|\mathcal{I} - \hat{\mathcal{I}}\|_{\mathcal{H}\to\mathbb{R}^d}}$$
$$= \alpha\|G^{-1}\|_{\mathbb{R}^d}^2 \|\mathcal{I} - \hat{\mathcal{I}}\|_{\mathcal{H}\to\mathbb{R}^d} + \mathcal{O}(\|\mathcal{I} - \hat{\mathcal{I}}\|_{\mathcal{H}\to\mathbb{R}^d}^2) \tag{S2}$$

whenever $\|\mathcal{I} - \hat{\mathcal{I}}\|_{\mathcal{H}\to\mathbb{R}^d}$ is sufficiently small.

Putting these results together we obtain the following upper bound for the Wasserstein distance:

$$W_2(\bar{\mu}, \hat{\mu}) \leq (C_{\mathcal{I},1} + C_{\mathcal{I},2})\|\mathcal{I} - \hat{\mathcal{I}}\|_{\mathcal{H}\to\mathbb{R}^d} + C_f\|f - \hat{f}\|_{\mathbb{R}^d} + \mathcal{O}\left(\|\mathcal{I} - \hat{\mathcal{I}}\|_{\mathcal{H}\to\mathbb{R}^d}^2\right)$$

$$C_{\mathcal{I},1} = \|C\|_{\mathcal{H}} \left[\left(\|a\|_{\mathcal{H}} + \|C\|_{\mathsf{HS}}^{\frac{1}{2}}\right)\left(\|G^{-1}\mathcal{I}\|_{\mathcal{H}\to\mathbb{R}^d} + \|\hat{G}^{-1}\hat{\mathcal{I}}\|_{\mathcal{H}\to\mathbb{R}^d}\right) + \|G^{-1}f\|_{\mathbb{R}^d}\right]$$

$$C_{\mathcal{I},2} = \alpha\|G^{-1}\|_{\mathbb{R}^d}^2 \|C\|_{\mathcal{H}} \left(\|G^{-1}f\|_{\mathbb{R}^d} + \|\mathcal{I}\|_{\mathcal{H}\to\mathbb{R}^d}\|\hat{\mathcal{I}}\|_{\mathcal{H}\to\mathbb{R}^d}\left(\|a\|_{\mathcal{H}} + \|C\|_{\mathsf{HS}}^{\frac{1}{2}}\right)\right)$$

$$C_f = \|C\|_{\mathcal{H}} \|G^{-1}\hat{\mathcal{I}}\|_{\mathcal{H}\to\mathbb{R}^d}.$$

$\square$

*Proof of Theorem 3.5.* We focus on the forward case, noting that the proof for the adjoint case is essentially identical. To prove the result we must relate the norms $\|\mathcal{I}_F - \hat{\mathcal{I}}_F\|_{\mathcal{H}\to\mathbb{R}^d}$ and $\|f_F - \hat{f}_F\|_{\mathbb{R}^d}$ to the assumed bounds on $\frac{\partial F}{\partial u}$ and $\frac{\partial F}{\partial p}$ from the statement of the corollary.

For the former, note that

$$\|\mathcal{I}_F - \hat{\mathcal{I}}_F\|_{\mathcal{H}\to\mathbb{R}^d} = \sup_{\substack{u \in \mathcal{H} \\ \|u\|=1}} \|(\mathcal{I}_F - \hat{\mathcal{I}}_F)u\|_2$$

$$= \sup_u \left[\sum_{j=1}^d ([\mathcal{I}_{F,j} - \hat{\mathcal{I}}_{F,j}]u)^2\right]^{\frac{1}{2}}$$

$$= \sup_u \left[\sum_{j=1}^d \left(\tilde{\mathcal{I}}_{F,j}\left[\frac{\partial F}{\partial u}[U(p_j), p_j] - \frac{\partial F}{\partial u}[\hat{U}(p_j), p_j]\right]u\right)^2\right]^{\frac{1}{2}}$$

$$\leq \left[\sum_{j=1}^d \left(\|\tilde{\mathcal{I}}_{F,j}\|_{\mathcal{F}\to\mathbb{R}}^2 \left\|\frac{\partial F}{\partial u}[U(p_j), p_j] - \frac{\partial F}{\partial u}[\hat{U}(p_j), p_j]\right\|_{\mathcal{U}\to\mathcal{F}}^2\right)\right]^{\frac{1}{2}}$$

$$\leq M\epsilon\sqrt{d}$$

Similarly, for the latter

$$\|f_F - \hat{f}_F\|_2 = \left[\sum_{j=1}^d \left(\tilde{\mathcal{I}}_{F,j}\left[\frac{\partial F}{\partial p}[U(p_j),p_j] - \frac{\partial F}{\partial p}[\hat{U}(p_j),p_j]\right]\right)\right]^{\frac{1}{2}}$$

$$\leq \left[\sum_{j=1}^d \|\tilde{\mathcal{I}}_{F,j}\|_{\mathcal{F}\to\mathbb{R}}^2 \left\|\frac{\partial F}{\partial p}[U(p_j),p_j] - \frac{\partial F}{\partial p}[\hat{U}(p_j),p_j]\right\|_{P\to\mathcal{F}}^2\right]^{\frac{1}{2}}$$

$$\leq M\epsilon\sqrt{d}.$$

Substituting these into the bound obtained in Theorem 3.4 yields the result. □

## S2 Posterior Updating

In this section we describe the procedure used to update the posterior distribution based on new information. For concreteness we will assume the notation of probabilistic forward sensitivity analysis; the approach followed for adjoint sensitivity analysis is essentially identical. Suppose we have a prior $\mu_F = \mathcal{N}(a_F, C_F)$ and have constructed the posterior $\bar{\mu}_F^1$ based upon the information $\mathcal{I}_F^1 = f_F^1$, where $\mathcal{I}_F^1 : \mathcal{U}_{\partial P} \to \mathbb{R}^d$ and $f_F^1 \in \mathbb{R}^d$. For computational purposes this requires the inverse of the matrix $\mathcal{I}_F^1 C(\mathcal{I}_F^1)^\dagger$ or, equivalently, the solution of linear systems involving this matrix. Thus, we suppose that the Cholesky factorisation $L_1$ has been computed, so that $\mathcal{I}_F^1 C(\mathcal{I}_F^1)^\dagger = L_1 L_1^\top$, where the factor $L_1$ is upper-triangular.

Now suppose that new information is supplied, $\mathcal{I}_F^2 : \mathcal{U}_{\partial P} \to \mathbb{R}^{d'}$ with corresponding information $f_F^2 \in \mathbb{R}^{d'}$, and suppose that $\mathcal{I}_F^2$ is linearly independent of $\mathcal{I}_F^1$. We wish to compute the posterior distribution $\bar{\mu}_F^2$ by conditioning the distribution $\bar{\mu}_F^1$ on this new information with as low a computational cost as possible.

A natural approach would be to simply perform the conditioning procedure described in Theorem 3.1 with $\mu_F = \bar{\mu}_F^1$, however in practise this approach was found to suffer from a high degree of numerical instability. Instead, we advocate conditioning the original prior on this new information by updating the Cholesky factorisation computed for $\mathcal{I}_F^1$. The form of the new Gramian matrix whose factorisation must be computed is

$$\mathcal{I}_F^{1:2} C(\mathcal{I}_F^{1:2})^\dagger = \begin{bmatrix} \mathcal{I}_F^1 C(\mathcal{I}_F^1)^\dagger & \mathcal{I}_F^1 C(\mathcal{I}_F^2)^\dagger \\ \mathcal{I}_F^2 C(\mathcal{I}_F^1)^\dagger & \mathcal{I}_F^2 C(\mathcal{I}_F^2)^\dagger \end{bmatrix}.$$

Following [Osborne, 2010, Appendix B], we may form the Cholesky factorisation of this matrix as

$$\mathcal{I}_F^{1:2} C(\mathcal{I}_F^{1:2})^\dagger = L_{1:2} L_{1:2}^\top$$

$$L_{1:2} = \begin{bmatrix} L_1 & S_{12} \\ 0 & S_{22} \end{bmatrix}$$

where $S_{12} = (L_1^\top)^{-1}\mathcal{I}_F^1 C(\mathcal{I}_F^2)^\dagger$, and $S_{22}$ is the Cholesky factorisation of a $\mathbb{R}^{d' \times d'}$ matrix given by
$$S_{22}S_{22}^\top = \mathcal{I}_F^2 C(\mathcal{I}_F^2)^\dagger - S_{12}S_{12}^\top.$$

Examining the cost of this procedure, we see that the only near-cubic operation required is this Cholesky factorisation, so that the update is approximately $\mathcal{O}((d')^3)$. This may then be used to compute the posterior from Theorem 3.1 without incurring the $\mathcal{O}((d+d')^3)$ cost that would be required to compute it *without* reusing the Cholesky factor. We note that a triangular solve of $\mathcal{O}((d+d')^2)$ is required to compute the full posterior distribution. Thus, as the size of the system grows, the cost will still increase at a quadratic rate. This was not found to be prohibitive for the experiments presented in this paper, and it is possible that with further optimisations this cost may also be reduced.

## S3 Gradient Descent Algorithm

Algorithm 1 provides pseudocode for the gradient descent algorithm described in Section 4.1.

---
**Algorithm 1** Pseudocode for the classical gradient descent procedure described in Section 4.1. Here $p^0$ is some starting point, $g$ and $\frac{dg}{dp}$ are functions returning the quantity to be minimised and its gradient, respectively, and $\epsilon$ is some demanded tolerance. The routine BLS is the backtracking line search method described in Nocedal and Wright [2006, Algorithm 3.1].

---
1: **procedure** GD$(p^0, g, \frac{dg}{dp}, \epsilon)$
2:     **for** $n = 1, 2, \ldots$ **do**
3:         $j^n \leftarrow \frac{dg}{dp}(p^{n-1})$
4:         $\gamma^n \leftarrow \text{BLS}(p^{n-1}, g, j^n)$
5:         **if** $\gamma^n < \epsilon$ **then**
6:             **return** $p^{n-1}$
7:         **end if**
8:         $p^n \leftarrow p^{n-1} - \gamma^n j^n$
9:     **end for**
10: **end procedure**

---

## S4 Additional Results for the Fitzhugh—Nagumo Example

The constraint function $F(u, p)$ is given by:
$$F(u, p) := \begin{bmatrix} \frac{du_1}{dt} - u_1 + \frac{u_1^3}{3} + u_2 + p_1 \\ \frac{du_2}{dt} + \frac{p_3 u_2 - u_1 - p_2}{p_4} \end{bmatrix}.$$

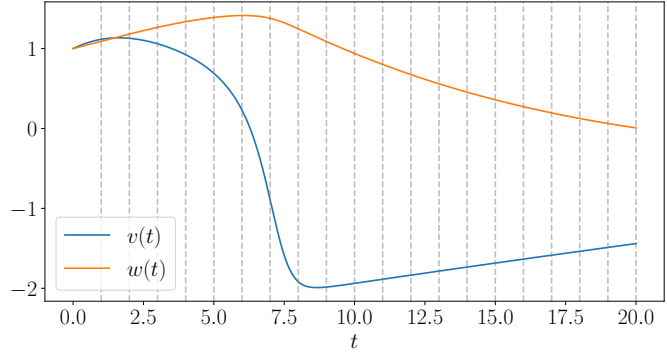

(a) Solutions $v, w$.

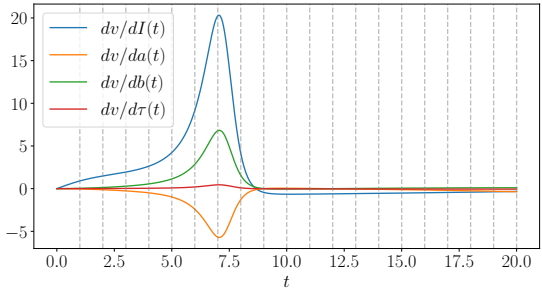

(b) $\frac{dv}{dp}$.

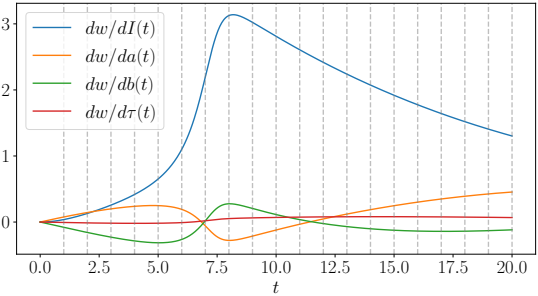

(c) $\frac{dw}{dp}$.

Figure S1: Top row: solutions $v, w$ to the Fitzhugh—Nagumo ODE as solved in Section 5.1, in the range $[0, 20]$. Bottom row: sensitivities of those solutions with-respect-to the four parameters $I, a, b, \tau$. The gray lines mark the points at which data was obtained.

Its derivatives $\frac{\partial F}{\partial u}$ and $\frac{\partial F}{\partial p}$ are:

$$\frac{\partial F}{\partial u}[u, p](u') = \begin{bmatrix} \frac{\partial u'_1}{\partial t} - u'_1 + u_1^2 u'_1 & u'_2 \\ -\frac{u'_1}{p_4} & \frac{\partial u'_2}{\partial t} + \frac{p_3 u'_2}{p_4} \end{bmatrix}$$

$$\frac{\partial F}{\partial p}[u, p] = \begin{bmatrix} 1 & 0 & 0 & 0 \\ 0 & -\frac{1}{p_4} & \frac{u_2}{p_4} & -\frac{p_3 u_2 - u_1 - p_2}{p_4^2} \end{bmatrix}$$

where the latter is simply a derivative in the classical sense, since $\dim(P) < \infty$, and we adopt the numerator layout convention for vector-by-vector derivatives.

Fig. S1 shows the solution and the sensitivities for the Fitzhugh—Nagumo ODE described in Section 5.1.

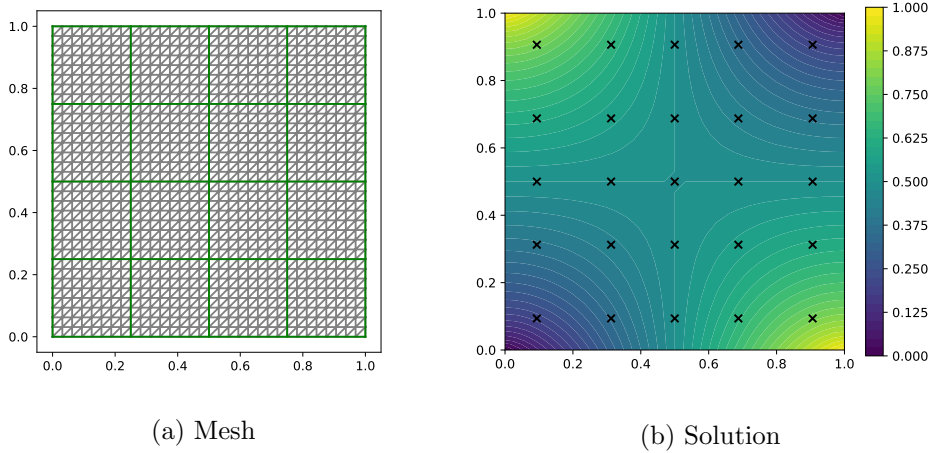

(a) Mesh  (b) Solution

Figure S2: Mesh and solution for the groundwater flow model from Section 5.2. In Fig. S2a the mesh used in the finite-element solver is depicted. The green lines bound the cells over which the piecewise constant parameter is defined; in this case the integer defining the parameter is $N = 4$, so there are 16 cells. In Fig. S2b the solution is shown, for the parameter $p = \mathbf{1}_{16}$. The points at which evaluations of the solution are taken to define the likelihood are shown as black crosses.

## S5 Additional Results for Groundwater Flow Example

Fig. S2 shows the mesh and a candidate solution for the groundwater flow problem described in Section 5.2.

## References

M. Gelbrich. On a formula for the L2 Wasserstein metric between measures on Euclidean and Hilbert spaces. *Mathematische Nachrichten*, 147(1):185–203, 1990. doi: 10.1002/mana.19901470121. URL https://doi.org/10.1002/mana.19901470121.

I. Gohberg, S. Goldberg, and M. A. Kaashoek. *Basic Classes of Linear Operators*. Birkhäuser Basel, 2003. doi: 10.1007/978-3-0348-7980-4.

V. Masarotto, V. M. Panaretos, and Y. Zemel. Procrustes metrics on covariance operators and optimal transportation of Gaussian processes. *Sankhya A*, 81(1):172–213, May 2018. doi: 10.1007/s13171-018-0130-1.

J. Nocedal and S. J. Wright. *Numerical Optimization*. Springer New York, 2006. doi: 10.1007/978-0-387-40065-5. URL https://doi.org/10.1007/978-0-387-40065-5.


\fi

\end{document}